\documentclass[10pt]{elsarticle}

\usepackage{printlen}
\usepackage[english]{babel}
\usepackage[utf8]{inputenc}
\usepackage{natbib}
\usepackage{hyperref}
\usepackage{verbatim}

\usepackage{graphicx}
\usepackage[caption=false]{subfig}
\usepackage{color}

\usepackage{amsmath,amsfonts,amssymb,amsthm}
\usepackage{mathtools}
\usepackage{empheq}
\usepackage{bm}
\usepackage[
  space-before-unit,
  range-units=repeat,
  detect-weight=true,
  detect-family=true
]{siunitx}

\usepackage{geometry}
\usepackage{pdflscape}

\usepackage{grffile}

\usepackage{multirow}
\usepackage{longtable}

\usepackage{algorithm}
\usepackage[]{algpseudocode}

\newcommand {\brc}[1]{\left( #1 \right)}

\newcommand {\txt}[1]{\text{#1}}
\newcommand {\abs}[1]{\left|{#1}\right|}

\newcommand {\nor}[1]{\left\lVert #1 \right\rVert}

\renewcommand {\vec}[1]{\boldsymbol #1 }

\graphicspath{{graphics/}}

\journal{ArXiv}

\algblock{ParFor}{EndParFor}
\algnewcommand\algorithmicparfor{\textbf{parfor}}
\algnewcommand\algorithmicpardo{\textbf{do}}
\algnewcommand\algorithmicendparfor{\textbf{end\ parfor}}
\algrenewtext{ParFor}[1]{\algorithmicparfor\ #1\ \algorithmicpardo}
\algrenewtext{EndParFor}{\algorithmicendparfor}

\algtext*{EndIf}%
\algtext*{EndFor}%
\algtext*{EndWhile}%
\algtext*{EndParFor}%

\begin{document}

\begin{frontmatter}


\title{Parallel Numerical Tensor Methods for High-Dimensional PDEs}

\author[stanford]{A. M. P. Boelens}
\ead{boelens@stanford.edu}
\author[ucsc]{D. Venturi}
\ead{venturi@ucsc.edu}
\author[stanford]{D. M. Tartakovsky}
\ead{tartakovsky@stanford.edu}

\address[stanford]{Department of Energy Resources Engineering, Stanford University, Palo Alto, CA 94305}
\address[ucsc]{Department of Applied Mathematics and Statistics, UC Santa Cruz, Santa Cruz, CA 95064}

\begin{abstract}
High-dimensional partial-differential equations (PDEs) arise in a 
number of fields of science and engineering, where they are used to describe the evolution of joint probability functions. Their examples include the Boltzmann and  Fokker-Planck 
equations. We develop new parallel algorithms to solve high-dimensional PDEs.
The algorithms are based on canonical and hierarchical numerical 
tensor methods combined with alternating least squares
and hierarchical singular value decomposition. Both implicit 
and explicit integration schemes are presented 
and discussed. We demonstrate the accuracy and 
efficiency of the proposed new algorithms in computing 
the numerical solution to both an advection 
equation in six variables plus time and a linearized version of the 
Boltzmann equation.
\end{abstract}

\begin{keyword}
\MSC[2010] 00-01\sep  99-00
\end{keyword}

\end{frontmatter}

\section{Introduction}

High-dimensional partial-differential equations (PDEs) arise in a 
number of fields of science and engineering. 
For example, they play an important role in modeling 
rarefied gas dynamics \cite{cercignani1988}, 
stochastic dynamical systems \cite{moss1995,sobczyk2001,venturi2012}, 
structural dynamics \cite{li2009,shlesinger1998},  
turbulence \cite{monin2007b,pope1994,frisch1995,venturi2018},  
biological networks \cite{cowen2017}, and 
quantum systems \cite{justin2002,amit2005}.
In the context of kinetic theory and stochastic dynamics, 
high-dimensional PDEs typically describe the evolution of 
a (joint) probability density function (PDF) of system states, 
providing macroscopic, continuum descriptions of Langevin-type stochastic 
systems. Classical examples of such PDEs include the Boltzmann 
equation \cite{cercignani1988} and the Fokker-Planck 
equation~\cite{risken1989}. More recently, PDF equations 
have been used to quantify uncertainty in model predictions 
\cite{tartakovsky2015}.  

High dimensionality and resulting computational complexity of 
PDF equations can be mitigated by using particle-based methods 
\cite{pope1994,dimarco2014}. 
Well known examples are direct simulation Monte-Carlo (DSMC) 
\cite{rjasanow2004} and the Nambu-Babovsky 
method \cite{Babovsky1986}. Such methods preserve physical 
properties of a system and exhibit high computational efficiency 
(they scale linearly with the number of particles), in particular in 
simulations far from statistical equilibrium \citep{ansumali2011,ramanathan2009}. 
Moreover, these methods have relatively low memory requirements 
since the particles tend to concentrate where the distribution function 
is not small. However, the accuracy of particle methods may 
be poor and their predictions are subject to significant statistical fluctuations 
\citep{ramanathan2009,peraud2014,rogier1994}. Such fluctuations need to be post-processed  appropriately, for 
example by using variance reduction techniques. Also, 
relevance and applicability of these numerical strategies to 
PDEs other than kinetic equations are not clear.
To overcome these difficulties, several general-purpose algorithms 
have been recently proposed to compute the numerical solution 
to rather general high-dimensional linear PDEs. The most efficient techniques, 
however, are problem specific \cite{cho2016,zhang2017,weinan2017,dimarco2014}. 
For example, in the context of kinetic methods, there is an extensive literature 
concerning the (six-dimensional) Boltzmann equation and its solutions close to 
statistical equilibrium. Deterministic methods for solving such problems 
include semi-Lagrangian schemes and discrete velocity models. 
The former employ a fixed computational grid, account for transport features of the 
Boltzmann equation in a fully Lagrangian framework, and usually adapt operator 
splitting. The latter employ a regular grid in velocity and a discrete collision
operator on the points of the grid that preserves the main physical properties
(see \S 3 and \S 4 in Di Marco and Pareschi \cite{dimarco2014} for an in-depth
review of semi-Lagrangian methods and discrete velocity models for kinetic
transport equations, respectively). 

We develop new parallel algorithms 
to solve high-dimensional partial differential equations
based on numerical tensor methods \cite{Bachmayr2016}. 
The algorithms we propose are based on canonical and 
hierarchical tensor expansions combined 
with alternating least squares and hierarchical 
singular value decomposition. The key element that opens 
the possibility to solve high-dimensional PDEs numerically with 
tensor methods is that tensor approximations proved to be 
capable of representing function-related $N$-dimensional data arrays 
of size $Q^N$ with log-volume complexity $\mathcal{O}(N\log(Q))$. 
Combined with traditional deterministic numerical 
schemes (e.g., spectral collocation method \cite{hesthaven2007}, 
these novel representations allow one to compute the 
solution to high-dimensional PDEs using low-parametric 
rank-structured tensor formats.

This paper is organized as follows. In Section \ref{sec:tensorseries} we
discuss tensor representations of $N$-dimensional functions. In particular, 
we discuss canonical tensor decomposition and hierarchical tensor networks, 
including hierarchical Tucker and tensor train expansions. We also address the 
problem of computing tensor expansion via alternating least squares. 
In Section \ref{sec:integrators} we develop several algorithms that rely on numerical 
tensor methods to solve high-dimensional PDEs. 
Specifically,  we discuss schemes based on implicit and explicit 
time integration, combined with dimensional splitting, 
alternating least squares and hierarchical singular value decomposition. 
In Section \ref{sec:numericalresults} we demonstrate the effectiveness 
and computational efficiency of the proposed new algorithms 
in solving both an advection equation in six 
variables plus time and a linearized version of the Boltzmann equation.

\section{Tensor Decomposition of High-Dimensional Functions}
\label{sec:tensorseries}

Consider a multivariate scalar field 
$f(\bm z): D\subseteq \mathbb R^N \rightarrow \mathbb R$. 
In this section we briefly review effective representations 
of $f$ based on tensor methods. 
In particular, we discuss the canonical tensor decomposition 
and hierarchical tensor methods, e.g., tensor train and 
hierarchical Tucker expansions.

\subsection{Canonical Tensor Decomposition}
\label{sec:CP}

The canonical tensor decomposition of the multivariate function
$f(\bm z): D\subseteq \mathbb R^N \rightarrow \mathbb R$ is a 
series expansion of the form 
\begin{equation}
f(\bm z)\simeq\sum_{l=1}^r \prod_{k=1}^NG^l_j(z_j),
\label{functional-SSE}
\end{equation}
where $G^l_i(z_i): \mathbb R \rightarrow \mathbb R$ are one-dimensional 
functions usually represented relative to a 
known basis $\{\phi_1,...,\phi_Q\}$, i.e., 
\begin{equation}
G^l_i(z_i)=\sum_{k=1}^Q\beta_{ik}^l \phi_k(z_i).
\label{gfun}
\end{equation}
The quantity $r$ in \eqref{functional-SSE} is 
called {separation rank}.
Although general/computable theorems relating a prescribed accuracy of the representation of $f(\bm z)$ to the value of the separation rank $r$ are still lacking, there are cases where the 
expansion \eqref{functional-SSE} is {exponentially 
more efficient} than one would expect 
a priori~\cite{beylkin2009}.


\paragraph{Alternating Least Squares (ALS) Formulation}
Development of robust and efficient algorithms to 
compute \eqref{functional-SSE} to any desired accuracy is 
still a relatively open question (see
\cite{acar2011,espig2012,karlsson2016,doostan2013,cho2016} 
for recent progresses). Computing the tensor 
components $G_k^l(z_k)$ usually relies on (greedy) 
optimization techniques, such as alternating least 
squares (ALS) \cite{reynolds2016,battaglino2017,acar2011,beylkin2009} 
or regularized Newton methods \cite{espig2012}, which are only locally 
convergent \cite{uschmajew2012}, i.e., the final result may depend 
on the initial condition of the algorithm. 

In the least-squares setting, the tensor 
components $G^l_j(z_j)$ are computed by minimizing a norm of the residual,
\begin{equation}
R(\bm z)=f(\bm z)-\sum_{l=1}^r \prod_{k=1}^N G^l_j(z_j), 
\label{Fresidual}
\end{equation}  
with respect to the $rNQ$ degrees of freedom $\beta_{hj}^n$
i.e., 
\begin{equation}
 \min_{\beta_{hj}^n}\left\|R(\bm z)\right\|.
 \label{als}
\end{equation}
Assuming $f(\bm z)$ to be periodic in the hyper-cube $D=[-b,b]^N$ ($b>0$),
 we define the norm $\|\cdot \|$ as a standard $L_2$ norm
\begin{align}
\left\| R(\bm z)\right\|_{L_2}^2 =& \int_{-b}^b\cdots \int_{-b}^b 
R(\bm z)^2 \text d\bm z,\nonumber\\
=& \int_{-b}^b\cdots \int_{-b}^b 
\left[f(\bm z) -\sum_{n=1}^r\prod_{k=1}^N \left(\sum_{s=1}^Q
\beta_{ks}^n\phi_s(z_k)\right)\right]^2   \text d\bm z.
\label{continuous_norm_W}
\end{align}
In the alternating least-squares (ALS) paradigm, 
we compute the minimizer of the residual \eqref{Fresidual}
by splitting the non-convex optimization problem \eqref{als} into 
a {\em sequence} of convex low-dimensional convex optimization
problems. To illustrate the method, let us define vectors 
\begin{align}
\bm \beta_1 = (\beta^1_{11}, ..., \beta^1_{1Q},...,\beta^r_{11}, ..., \beta^r_{1Q})^T, \quad\cdots,\quad 
\bm \beta_N=(\beta^1_{N1}, ..., \beta^1_{NQ},...,\beta^r_{N1}, ..., \beta^r_{NQ})^T. 
\end{align}
Each $\bm \beta_i$ collects the degrees of freedom 
of all functions $\{G^1_{i}(z_i), ..., G^r_i(z_i)\}$ depending 
on $z_i$. Next, the optimization 
problem \eqref{als} is split into a sequence of 
convex optimization problems,
\begin{equation}
\min_{\bm\beta_1}\left\| R\right\|_{L_2}^2, 
\quad\cdots, \quad 
\min_{\bm\beta_N}\left\| R\right\|_{L_2}^2.
\label{als11}
\end{equation}
This sequence is not equivalent to the full problem \eqref{als} because, 
in general, it does not allow one to compute the global minimizer of \eqref{als} \cite{uschmajew2012,espig2015,bezdek2003,rohwedder2013}. 
The Euler-Lagrange equations associated 
with \eqref{als11} are of the form\footnote{
Recall that minimizing the residual \eqref{Fresidual}
with respect to $\beta_{hj}^n$ 
is equivalent to impose orthogonality 
relative to the space spanned by the functions
\begin{equation}
\phi_h(z_j)\prod^N_{\substack{k=1\\k\neq j}}G^{n}_k(z_k).
\end{equation}
}
\begin{equation}
\bm A_j \bm \beta_j =\bm g_j\qquad j=1,\ldots,N, 
\label{SSYSTEM}
\end{equation}
where
\begin{equation}
\bm A_j=\left[
\begin{array}{cccc}
A^{11}_{j11}\cdots A^{11}_{j1Q} & A^{12}_{j11}\cdots A^{12}_{j1Q} 
& \cdots &A^{1r}_{j11}\cdots A^{1r}_{j1Q}\\
\vdots & \vdots & \vdots & \vdots \\
A^{11}_{jQ1}\cdots A^{11}_{jQQ} & A^{12}_{jQ1}\cdots A^{12}_{jQQ} 
&\cdots & A^{1r}_{jQ1}\cdots A^{1r}_{jQQ}\\
\vdots & \vdots & \vdots & \vdots \\
A^{r1}_{j11}\cdots A^{r1}_{j1Q} & A^{r2}_{j11}\cdots A^{r2}_{j1Q} 
&\cdots & A^{rr}_{j11}\cdots A^{rr}_{j1Q}\\
\vdots & \vdots & \vdots & \vdots \\
A^{r1}_{jQ1}\cdots A^{r1}_{jQQ} & A^{r2}_{jQ1}\cdots A^{r2}_{jQQ} 
&\cdots & A^{rr}_{jQ1}\cdots A^{rr}_{jQQ}
\end{array}
\right],\qquad 
\bm g_j=\left[
\begin{array}{c}
 g_{j1}^1\\
 \vdots\\
g_{jQ}^1\\
 \vdots\\
 g_{j1}^r\\
 \vdots\\
g_{jQ}^r\\ 
\end{array}
\right],
\label{ORDERING}
\end{equation}
and 
\begin{equation}
A^{ln}_{jhs}=\int_{-b}^b \phi_h(z_j) \phi_s(z_j) \text dz_j 
\prod_{\substack{k=1\\k\neq j}}^N
\int_{-b}^b G_k^l(z_k)G_k^n(z_k) \text dz_k,
\label{A^{ln}_{jhp}}
\end{equation}
\begin{equation}
g_{jh}^n = \int_{-b}^b\cdots \int_{-b}^b f(\bm x) \phi_h(x_j)\prod_{\substack{k=1\\k\neq j}}^N
G_k^n(x_k) \text dz_1\cdots \text dz_N.
\label{f_{jh}^n}
\end{equation}
The matrices $\bm A_j$ are symmetric, positive 
definite and of size $rQ\times rQ$.

\paragraph{Convergence of the ALS Algorithm} 
The ALS algorithm described above is  
an alternating optimization scheme, i.e.,  
a nonlinear block Gauss--Seidel 
method (\cite{ortega1970}, \S 7.4). 
There is a well--developed local 
convergence theory for this type of methods~\cite{ortega1970,bezdek2003}. 
In particular, it can be shown that ALS is locally 
equivalent to the linear block Gauss--Seidel 
iteration applied to the Hessian matrix. 
This implies that ALS is linearly convergent in the iteration 
number \cite{uschmajew2012}, 
provided that the Hessian of the residual 
is positive definite (except on a trivial null 
space associated with the scaling non-uniqueness 
of the canonical tensor decomposition). The last assumption 
may not be always satisfied. Therefore, convergence 
of the ALS algorithm cannot be granted in general.
Another potential issue of the ALS algorithm is 
the poor conditioning of the matrices $\bm A_j$ 
in \eqref{SSYSTEM}, which can addressed by 
regularization \cite{reynolds2016,battaglino2017}.
The canonical tensor decomposition 
\eqref{functional-SSE} in $N$ dimensions has relatively 
small memory requirements. In fact, the number of degrees of freedom 
that we need to store is $rN  Q$, where $r$ is 
the separation  rank, and $Q$ is the number of degrees of 
freedom employed in each tensor component $G_k^l(z_k)$.
Despite the relatively low-memory requirements, it is often desirable 
to employ scalable parallel versions the ALS algorithm 
\cite{karlsson2016,kaya2018} to compute the canonical 
tensor expansion \eqref{functional-SSE} 
(see Appendix \ref{app:parallelization})

\subsection{Hierarchical Tensor Methods}
\label{sec:Htensor}
\begin{figure}[ht]
\centerline{\hspace{-1.cm}\footnotesize Hierarchical Tucker\hspace{4.2cm} Tensor Train}
\centerline{
\includegraphics[width=6.5cm]{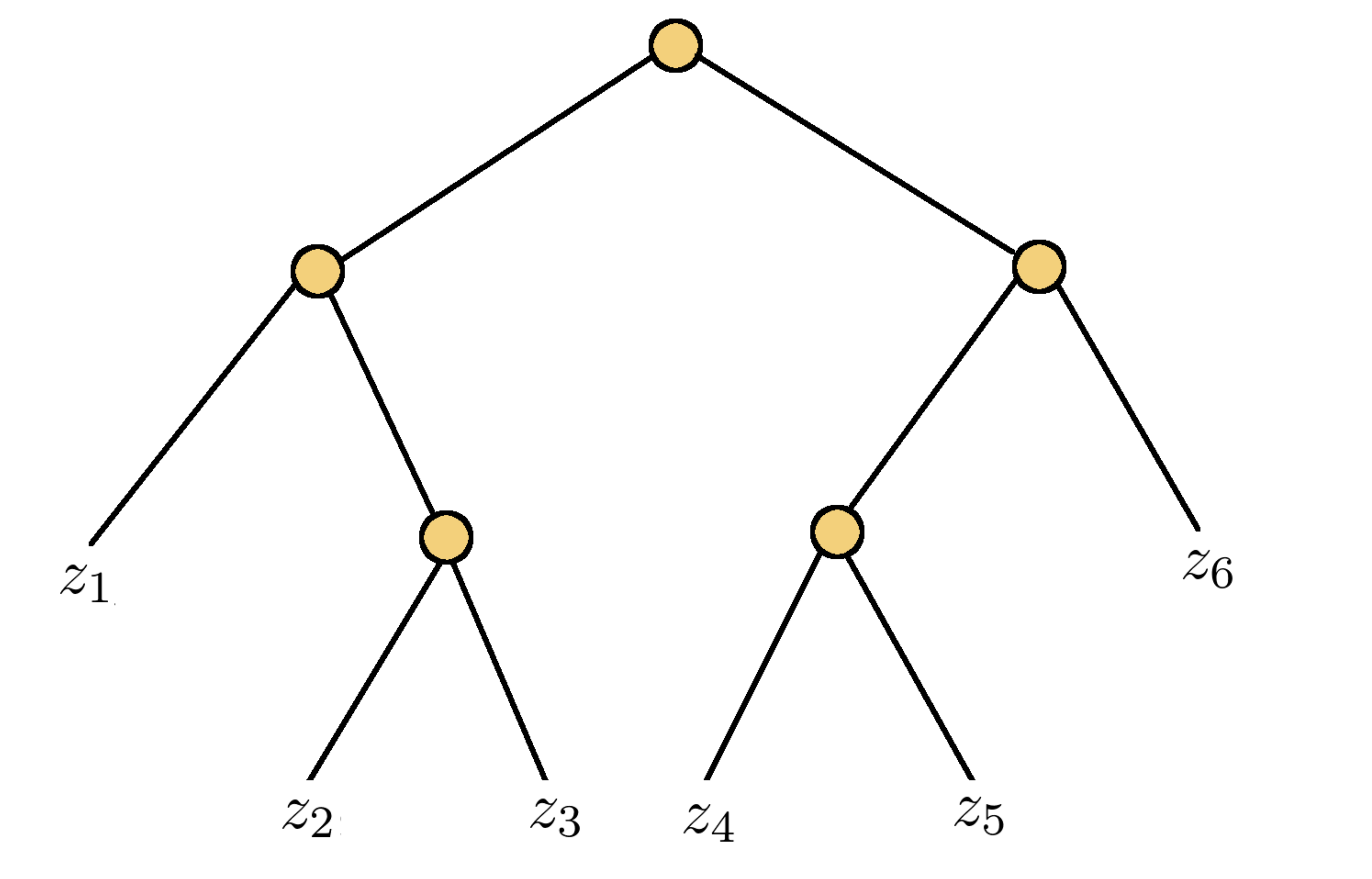}
\hspace{1cm}
\includegraphics[width=6cm]{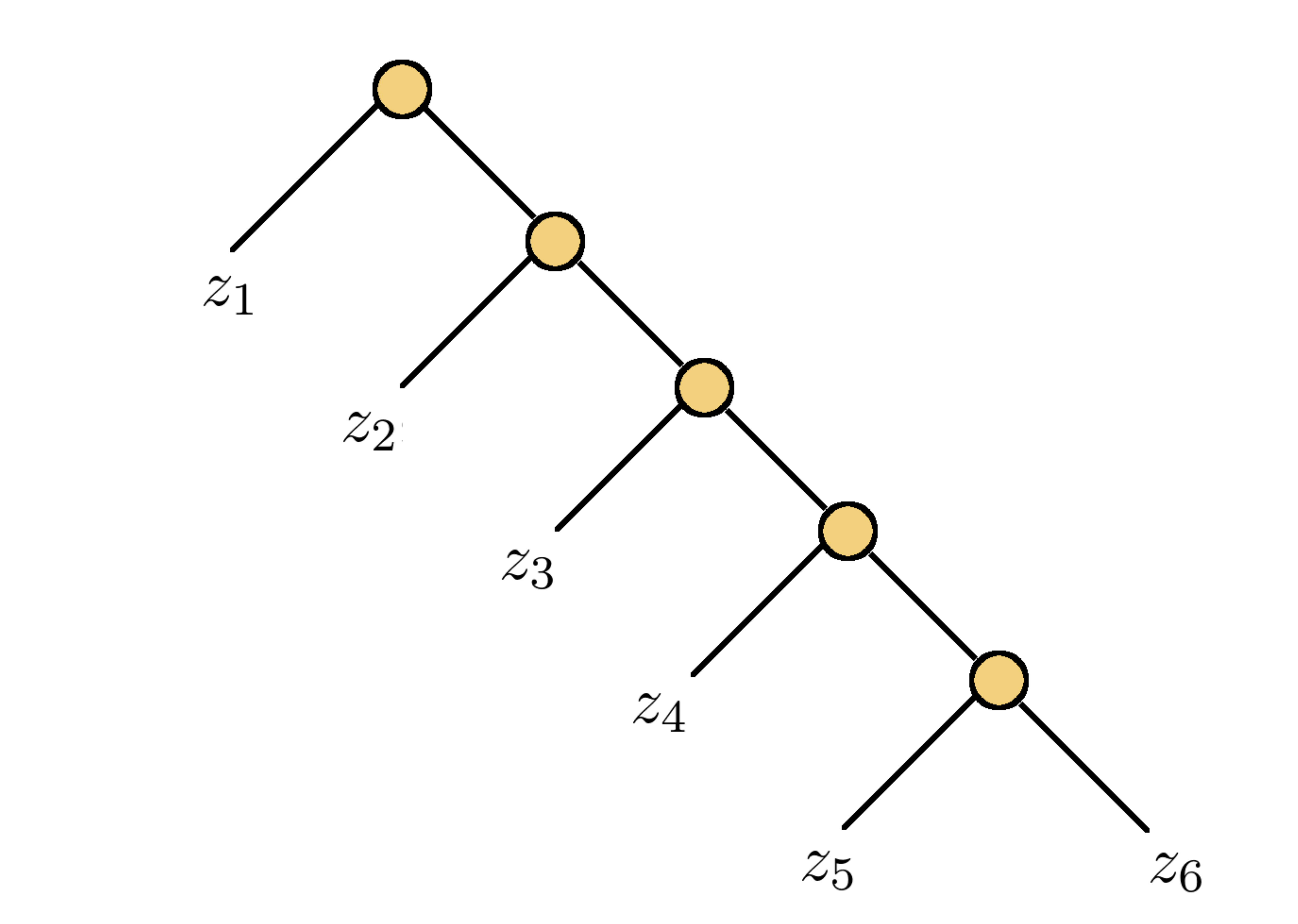}
}
\caption{Graph representation of the hierarchical Tucker (HT) and tensor train (TT) decomposition of a six-dimensional function $f(z_1,\ldots,z_6)$.}
\label{fig:HT_TT}
\end{figure}

Hierarchical Tensor methods \cite{Hackbusch2009,hackbusch2012}  
were introduced to mitigate the dimensionality problem in the  
core tensor of the classical Tucker decomposition \cite{kolda2009}. 
A key idea is to perform a sequence of 
Schmidt decompositions (or multivariate SVDs \cite{grasedyck2010,lathauwer2000}) until the 
approximation problem is reduced to a product 
of one-dimensional functions/vectors. 
To illustrate the method in a simple way, 
consider a six-dimensional function
$f(\bm z)= f(z_1,\ldots,z_6)$. 
We first split the variables as $(z_1,z_2,z_3)$ and 
$(z_4,z_5,z_6)$ through one Schmidt 
decomposition \cite{venturi2011} as
\begin{align}
f(\bm z)= &\sum_{i_7,i_8=1}^r A^{\{1\}}_{i_7i_8} T^{\{1,2,3\}}_{i_7}
(z_1,z_2,z_3)T^{\{4,5,6\}}_{i_8}(z_4,z_5,z_6).
\end{align}
Then we decompose $T^{\{1,2,3\}}_{i_7}(z_1,z_2,z_3)$ and 
$T^{\{4,5,6\}}_{i_8}(z_4,z_5,z_6)$ further by 
additional Schmidt expansions to obtain
 \begin{align}
f(\bm z)=&\sum_{i_7,i_8=1}^r A^{\{1\}}_{i_7i_8} 
\sum_{i_1,i_9=1}^r  A^{\{2\}}_{i_7i_1i_9}
 T^{\{1\}}_{i_1}(z_1)
 T^{\{2,3\}}_{i_9}(z_2,z_3)
\sum_{i_4,i_{10}=1}^r A^{\{3\}}_{i_8 i_4i_{10}}
T^{\{4\}}_{i_4}(z_4)T^{\{5,6\}}_{i_{10}}(z_5,z_6),\nonumber\\
&=\sum_{i_1,\cdots, i_{6}=1}^r
C_{i_1\cdots i_6}
T^{\{1\}}_{i_1}(z_1)T^{\{2\}}_{i_2}(z_2)T^{\{3\}}_{i_3}(z_3)
T^{\{4\}}_{i_4}(z_4)T^{\{5\}}_{i_5}(z_5)T^{\{6\}}_{i_{6}}(z_6).
\label{HTD}
\end{align}
The $6$-dimensional {\em core tensor}\footnote{A diagonalization 
of the core tensor in \eqref{HTD} would minimize the number of 
terms in the series expansion. Unfortunately, this is impossible for a tensor with 
dimension larger than 2 (see, e.g.,~\cite{peres1995,silva2008,kolda2009,hillar2013})
and, for complex tensors,~\cite{vannieuwenhoven2014}). A closer look at the canonical tensor 
decomposition \eqref{functional-SSE} reveals that such an expansion is 
in the form of a fully diagonal high-order Schmidt decomposition, i.e., 
\begin{equation}
f(x_1,\ldots,x_N)=\sum_{l=1}^{r} C_{l\cdots l}
G^{l}_{1}(x_1)\cdots G^{l}_N(x_N).
\end{equation}
The fact that diagonalization of $C_{l_1\cdots  l_6}$ is 
impossible in dimension larger than 2 implies that it is impossible 
to compute the canonical tensor decomposition of $f$ by standard 
linear algebra techniques.} 
is explicitly obtained as
\begin{equation}
C_{i_1\cdots i_{6}} =
\sum_{i_7,i_8=1}^r A^{\{1\}}_{i_7i_8} 
\sum_{i_9=1}^r  A^{\{2\}}_{i_7i_1i_9} A^{\{4\}}_{i_9i_2i_3}
\sum_{i_{10}=1}^r A^{\{3\}}_{i_8 i_4i_{10}}A^{\{5\}}_{i_{10}i_5i_6}.
\label{coreT}
\end{equation}
The procedure just described forms the foundation of 
hierarchical tensor methods. The key element is that 
the core tensor $C_{i_1\cdots i_{6}}$ resulting from 
this procedure is always factored as 
a {\em product of at most three-dimensional matrices}. 
This is also true in arbitrary dimensions. 
The tensor components $T^{\{k\}}_{i_j}$ and the factors of the  
core tensor \eqref{coreT} can be effectively 
computed by employing hierarchical singular value decomposition  \cite{grasedyck2010,Kressner2014,lathauwer2000}.
Alternatively, one can use an optimization framework 
that leverages recursive subspace factorizations 
\cite{DaSilva2015}. 
If we single out one variable at the time and perform a 
sequential Schmidt decomposition of the remaining 
variables we obtain the so-called {\em tensor-train} (TT)
decomposition \cite{oseledets2011,rohwedder2013}. 
Both tensor train and hierarchical tensor expansions
can be conveniently visualized by  
{\em graphs} (see Figure \ref{fig:HT_TT}). This is done 
by adopting the following standard rules: i)  
a node in a graph represents a tensor in as many variables as 
the number of the edges connected to it, ii) 
connecting two tensors by an edge represents 
a tensor contraction over the index 
associated with a certain variable. 

Efficient algorithms to perform basic operations 
between hierarchical tensors, such as addition, 
orthogonalization, rank reduction, scalar products, 
multiplication, and linear transformations, are discussed in 
\cite{grasedyck2015,kolda2009,Kressner2014,Nouy2017}.
Parallel implementations of such algorithms were recently proposed in 
\cite{etter2016,grasedyck2017}.
Methods for reducing the computational cost of tensor 
trains are discussed in \cite{Phan2013b,Zhou2013,Nouy2017}. 
Applications to the Vlasov kinetic equation can be found in 
\cite{Hatch2012,Kormann2015,Dolgov2014}.

\section{Numerical Approximation of  High-Dimensional PDEs}
\label{sec:integrators}

In this section we develop efficient numerical methods to integrate 
high-dimensional linear PDEs of the form 
\begin{equation}
\frac{\partial f(\bm z,t)}{\partial t} = L(\bm z) f(\bm z,t),
\label{PDE0}
\end{equation}
where $L(\bm z)$ is a linear operator. 
The algorithms we propose are based on numerical tensor 
methods and appropriate rank-reduction techniques, such as 
hierarchical singular value decomposition~\cite{grasedyck2010,grasedyck2017} and 
alternating least squares~\cite{karlsson2016,espig2012}.
To begin with, we assume that $L(\bm z)$ is a separable 
linear operator with separation rank $r_L$, i.e.,  an operator of the form
\begin{equation}
L(\bm z) = \sum_{l=1}^{r_L}\alpha_l L_1^l(z_1)\cdots L_N^l(z_N).
\label{Lz}
\end{equation}
For each $i$ and $j$, $L_i^j(z_i)$ is a linear operator acting only on variable $z_i$.
As an example,  the operator
\begin{equation}
L(z_1,z_2)= z^2_1\frac{\partial }{\partial z_1}-\sin(z_2)\frac{\partial^2 }{\partial z_1\partial z_2}.
\end{equation}
is separable with separation rank $r_L=2$ in dimension $N=2$. It can be written in the form \eqref{Lz} if we set  
$\alpha_1=\alpha_2=1$ and 
\begin{equation}
L_1^1(z_1)=z_1^2\frac{\partial }{\partial z_1}, \quad L_2^1(z_2)=1, \quad
L_1^2(z_1)=\frac{\partial }{\partial z_1}, \quad L_2^2(z_2)=-\sin(z_2)\frac{\partial }{\partial z_2}.
\end{equation}

\subsection{Tensor Methods with Implicit Time Stepping}
\label{sec:BoltzmanImplicit}
\begin{figure}[t]
\centerline{\includegraphics[height=5.2cm]{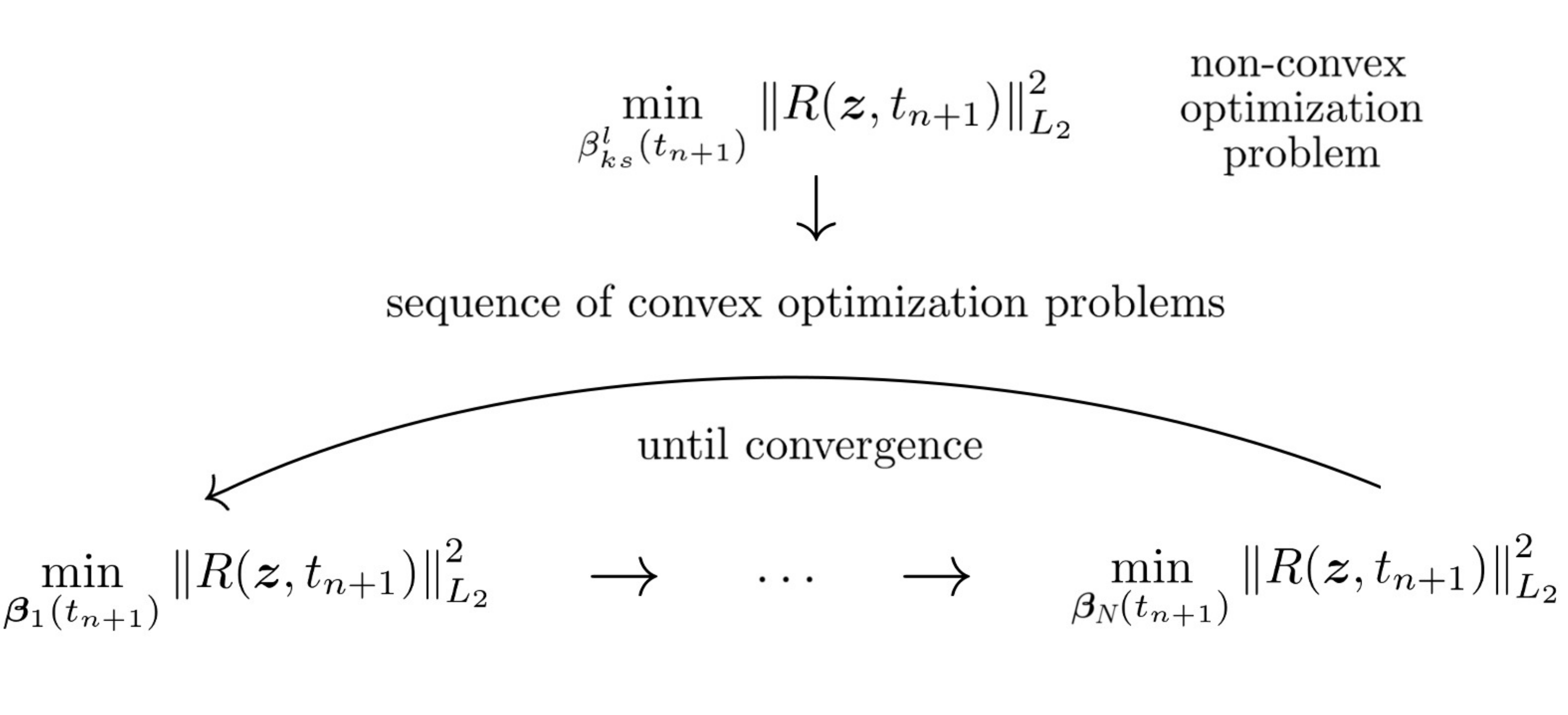}}
\caption{Sketch of the alternating least squares (ALS) algorithm for 
the minimization of the norm of the residual $R(\bm z,t_{n+1})$
defined in equation \eqref{rRr1}.}
\label{fig:ALS}
\end{figure}

Let us discretize the PDE in~\eqref{PDE0} in time by using the 
Crank-Nicolson method. To this end, consider an evenly-spaced 
grid $t_n= n\Delta t$ ($n=0,1,...$) with time step $\Delta t$.  This yields 
\begin{equation}
\left[I-\frac{\Delta t}{2} L(\bm z) \right]f_{n+1}(\bm z)=
\left[I+\frac{\Delta t}{2} L(\bm z) \right]f_n (\bm z)+ \Delta t\tau_{n+1}(\bm z),
\label{CNFDE}
\end{equation}
where $I$ is the identity matrix, $f_n(\bm z)=f(\bm z,t_n)$, and $\tau_{n+1}(\bm z)$ is the 
local truncation error of the Crank-Nicolson method 
at time $t_{n+1}$ \cite{quarteroni2007}.
Rewriting this in a more compact notation yields~\cite{cho2016}
\begin{equation}
A(\bm z) f_{n+1}(\bm z)=B(\bm z) f_n(\bm z)+\Delta t \tau_{n+1}(\bm z),
\end{equation}
where  
\begin{equation}
A = I-\frac{\Delta t}{2} L(\bm z), \quad \textrm{and}\quad 
B = I+\frac{\Delta t}{2} L(\bm z).
\end{equation}
It follows from the definition of $L(\bm z)$ in \eqref{Lz} that both $A(\bm z)$ and $B(\bm z)$ are separable operators 
of the form  
\begin{align}
A= \sum_{q=0}^{r_L} \eta_q E^q_1(z_1)\cdots E^q_N(z_N), \qquad 
B= \sum_{q=0}^{r_L} \zeta_q E^q_1(z_1)\cdots E^q_N(a_N),\label{BLs}
\end{align}
where $\eta_0=\zeta_0=1$, $E^0_j(z_j)=1$ ($j=1,...,N$), 
\begin{equation}
\eta_j=-\frac{\Delta t}{2},\qquad \zeta_j=-\eta_j,\qquad E^q_j(z_j)=L^q_j(z_j), \qquad q=1,...,r_L, \quad j=1,...,N. 
\end{equation}
A substitution of the canonical tensor  
decomposition\footnote{Recall that the 
functions $G_{k}^l(z_k,t_n)$ 
are in the form 
\begin{equation}
G_{k}^l(z_k,t_n)=\sum_{s=1}^Q 
\beta^l_{ks}(t_n)\phi_s(z_k),
\end{equation}
where $\beta^l_{ks}(t_n)$ ($l=1,...,r$, $k=1,...,N$, $s=1,...,Q$) 
are the degrees of freedom. }
\begin{equation}
\hat{f}_{n+1}(\bm z) = \sum_{l=1}^r \prod_{k=1}^N G^l_k(z_k,t_{n+1})
\label{CP-numerical0}
\end{equation}
into~\eqref{CNFDE} yields the residual 
\begin{equation}
R(\bm z,t_{n+1}) = A(\bm z)\hat{f}_{n+1}(\bm z)-B(\bm z)\hat{f}_n(\bm z),
\label{rRr1}
\end{equation}
in which the local truncation error $\Delta t \tau_{n+1}$ is embedded into $R(\bm z,t_{n+1})$.
Next, we minimize the residual using the alternating least squares
algorithm, described in Section \ref{sec:BoltzmanImplicit}, to obtain the solution at time $t_{n+1}$.
In particular, we look for a minimizer of \eqref{rRr1} 
computed in a {parsimonious way}. 
The key idea again is to split the optimization problem
\begin{equation}
\min_{\beta_{ks}^l(t_{n+1})} \left\|R(\bm z,t_{n+1})\right\|_{L_2}^2
\label{resnorm}
\end{equation}
into a sequence of optimization problems of smaller
dimension, which are solved sequentially and in parallel 
\cite{karlsson2016} (see Appendix \ref{app:parallelization}). 
To this end, we define
\begin{equation}
\bm \beta_k(t_n)=[\beta^1_{k1}(t_n), ..., \beta^1_{kQ}(t_n),...,
\beta^r_{k1}(t_n), ..., \beta^r_{kQ}(t_n)]^T \qquad k=1,...,N.
\end{equation}
The vector $\bm \beta_k(t_n)$ collects 
the degrees of freedom representing the solution 
functional along $z_k$ at time $t_n$, i.e., 
the set of functions $\{G_k^1(z_k,t_n),...,G^r_k(z_k,t_n)\}$.
Minimization of \eqref{resnorm} with respect to 
{ independent} variations of $\bm \beta_k(t_{n+1})$
yields a sequence of {\em convex} optimization problems (see Figure~\ref{fig:ALS})
\begin{equation}
\min_{\bm \beta_1(t_{n+1})} \left\|R(\bm z,t_{n+1})\right\|_{L_2}^2, \quad 
\min_{\bm \beta_2(t_{n+1})} \left\|R(\bm z,t_{n+1})\right\|_{L_2}^2, \quad 
\cdots,\quad 
\min_{\bm \beta_N(t_{n+1})} \left\|R(\bm z,t_{n+1})\right\|_{L_2}^2.
\label{sequential_min}
\end{equation}
\noindent
This set of equations defines the alternating least-squares (ALS) method. 
The Euler-Lagrange equations, which identify stationary points 
of \eqref{sequential_min}, are linear systems in the form
\begin{equation}
 \bm M^L_q \bm \beta_q(t_{n+1}) = \bm M^R_q \bm \beta_q(t_{n})\qquad q=1,...,N
\label{ALS2}
\end{equation}
where
\begin{equation}
\bm M^L_q = \sum_{e,z=0}^{r_L} K^L_{ez}
\left[
{\displaystyle \substack{\vspace{0.2cm}\\\vspace{0.04cm} N\\\bigcirc\\k=1\\k\neq q}}
\hat{\bm \beta}_{k}(t_{n+1})^T \bm 
E^{ez}_k \hat {\bm \beta}_{k}(t_{n+1})
\right]^T\otimes \left[{\bm E}^{ez}_q\right]^T, 
\label{ML}
\end{equation}
\begin{equation}
\bm M^R_q = \sum_{e,z=0}^{r_L} {K}^R_{ez}
\left[ 
{\substack{\vspace{0.2cm}\\\vspace{0.04cm} N\\\bigcirc\\k=1\\k\neq q}}
\hat{\bm \beta}_{k}(t_{n})^T 
{\bm E}^{ez}_k \hat{\bm \beta}_{k}(t_{n+1})
\right]^T\otimes \left[{\bm E}^{ez}_q\right]^T, 
\label{MR}
\end{equation}
and
\begin{align}
\left[{\bm E}^{ez}_q\right]_{sh}=\int_{-b}^{b} E_{q}^e(a)\phi_s(a) E_{q}^z(a)\phi_h(a) \text da.
\label{eq:77}
\end{align}
In~\eqref{ML}--\eqref{MR}, 
$\bigcirc$ denotes the Hadamard matrix 
product, $\otimes$ is the Kronecker matrix product, 
$\hat {\bm \beta}_{k}$ is the matricization of $\bm \beta_k$, i.e., 
\begin{equation}
\hat {\bm \beta}_{k}(t_n)=
\left[\begin{array}{ccc}
\beta^1_{k1}(t_n)&\cdots& \beta^r_{k1}(t_n)\\
\vdots  &\ddots & \vdots \\
\beta^1_{kQ}(t_n) &\cdots & \beta^r_{kQ} (t_n)
\end{array}
\right],
\end{equation}
$\bm E^{ez}_q$ are $Q\times Q$ matrices, 
and  $K^L_{ez}$ and ${K}^R_{ez}$ are entries of the 
matrices
\begin{equation}
\bm K^L = \bm \eta \bm \eta^T,\quad \textrm{and}\quad 
\bm K^R = \bm \eta \bm \zeta^T,
\end{equation}
where $\bm \eta$ and $\bm \zeta$ are column vectors with entries $\eta_q$ and $\zeta_q$ defined in~\eqref{BLs}.
The ALS algorithm, combined with 
canonical tensor representations, effectively reduces evaluation of $N$-dimensional integrals to evaluation of the sum of products of one-dimensional integrals.

\paragraph{Summary of the Algorithm}
Often many of the integrals in \eqref{eq:77} only differ by a prefactor and in that
case the number of unique integrals can be reduced to a very small number
($4$ in the case of the Boltzmann-BGK
equation discussed in Section \ref{sec:BoltzmanDiscretization}). 
Thus, it is convenient to precompute a map between such 
integrals and any entry in \eqref{eq:77}.
Such a map allows us to rapidly compute each 
matrix $\bm E_q^{ez}$  
in \eqref{ML} and \eqref{MR}, and therefore to efficiently build the matrix system. 
Next, we compute the canonical tensor 
decomposition of the initial condition $f(\bm z,0)$
by applying the methods described in Section \ref{sec:CP}.
This yields the set of vectors 
$\{\bm \beta_1(t_0),...,\bm \beta_N(t_0)\}$, from which the matrices 
$\bm M^L_1$ and $\bm M^R_1$ in \eqref{ML} and 
\eqref{MR} are built. To this end, we need an initial guess 
for $\{\bm \beta_1(t_1),...,\bm \beta_N(t_1)\}$ given, e.g., by $\{\bm \beta_1(t_0),...,\bm \beta_N(t_0)\}$ or its small random perturbation. 
With $\bm M^L_1$, $\bm M^R_1$ in place, 
we solve the linear system in~\eqref{ALS2} and 
update $\bm \beta_1(t_1)$. Finally, we recompute  $\bm M^L_2$ and $\bm M^R_2$  
(with the updated $\bm \beta_1(t_1)$) and solve for $\bm \beta_2(t_1)$. 
This process is repeated for $q=3,\ldots,N$ and iterated among 
all the variables until convergence (see Figure \ref{fig:ALS}). 
The ALS iterations are said to converge if 
\begin{equation}
\left\|\bm \beta^{(j)}_k(t_n)-\bm \beta^{(j-1)}_k(t_n)\right\|_2^2\leq 
\epsilon_{\textrm{Tol}}, \qquad j=1,\ldots,N_\text{max},
\end{equation}
where $j$ here denotes the iteration number whose maximum value is $N_\text{max}$, and $\epsilon_{\textrm{Tol}}$ is a prescribed tolerance on the increment of such iterations. 
Parallel versions of the ALS algorithm associated with canonical 
tensor decompositions were recently proposed in \cite{karlsson2016,kaya2018} 
(see also Appendix \ref{app:parallelization}).
We emphasize that the use of a hierarchical tensor series 
instead of \eqref{CP-numerical0} (see Section \eqref{sec:Htensor}) also allows one to compute the unknown degrees of the expansion by using parallel ALS algorithms. This question was recently addressed in~\cite{etter2016}.

\subsection{Tensor Methods with Explicit Time Stepping}
\label{sec:BoltzmanExplicit}
Let us discretize PDE \eqref{PDE0} in time by using any explicit 
time stepping scheme, for example the second-order Adams-Bashforth 
scheme \cite{quarteroni2007}
\begin{equation}
f_{n+2}(\bm z)=f_{n+1}(\bm z)+\frac{\Delta t}{2} L(\bm z)\left(3f_{n+1}(\bm z)-f_n(\bm z)\right) +\Delta t \tau_{n+2}(\bm z).
\label{AB2_}
\end{equation}
In this setting, the only operations needed to 
compute $f_{n+2}$ with tensor methods are: i) addition, 
ii) application of a (separable) linear operator $L$, and iii) rank 
reduction\footnote{On the other hand, the use of
implicit time-discretization schemes requires one to  
 solve linear systems with tensor methods~\cite{grasedyck2017,Kressner2014,reynolds2016}.}, 
the last operation being the most important among the three
\cite{kolda2009,Kressner2014,reynolds2016}. 
Rank reduction can be achieved, e.g., by hierarchical singular value 
decomposition \cite{grasedyck2010,grasedyck2017} or by 
alternating least squares \cite{karlsson2016}.

From a computational perspective, it is useful to 
split the sequence of tensor operations yielding $f_{n+2}$ in \eqref{AB2_} 
into simple tensor operations followed by rank reduction. For example, 
 $L(3f_{n+1}-f_n)$ might be evaluated as follows: i) compute 
$w_{n+1}=f_{n+1}-f_n$; ii) perform rank reduction on $w_{n+1}$; iii) compute 
$L w_{n+1}$ and perform rank reduction again.  
This allows one to minimize the memory requirements and overall 
computational cost.\footnote{Recall that adding two tensors 
with rank $r_1$ and $r_2$ usually yields a tensor of rank $r_1+r_2$.}
Unfortunately, splitting tensor operations into sequences of tensor operations 
followed by rank reduction can produce severe cancellation 
errors. In some cases this problem can 
be mitigated. For example, an efficient and robust algorithm, which 
allows us to split sums and rank reduction operations, was recently proposed 
in \cite{Kressner2014} in the context of 
hierarchical Tucker formats. The algorithm leverages  
the block diagonal structure arising from addition of
hierarchical Tucker tensors.

\section{Numerical Results}
\label{sec:numericalresults}
In this section, we demonstrate the accuracy and 
effectiveness of the numerical tensor methods discussed above. Two case studies are considered: 
an advection equation in six dimensions plus time and the Boltzmann-BGK equation.

\subsection{Advection Equation}
\begin{figure}[t]
\centerline{\footnotesize$t=0$\hspace{4.5cm}$t=3$\hspace{4.5cm}$t=6$}
\centerline{
\includegraphics[height=4.2cm]{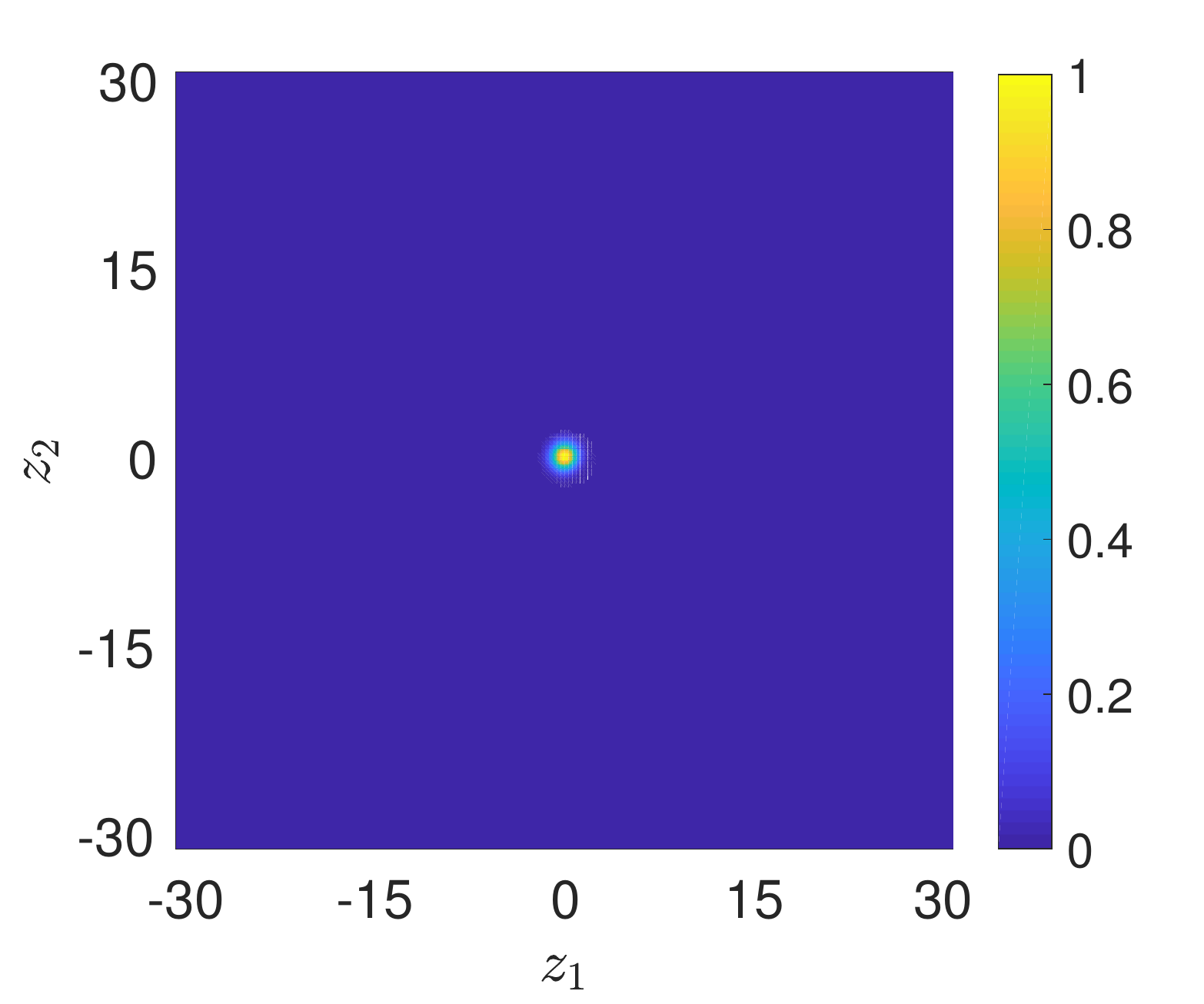}\hspace{-0.0cm}
\includegraphics[height=4.2cm]{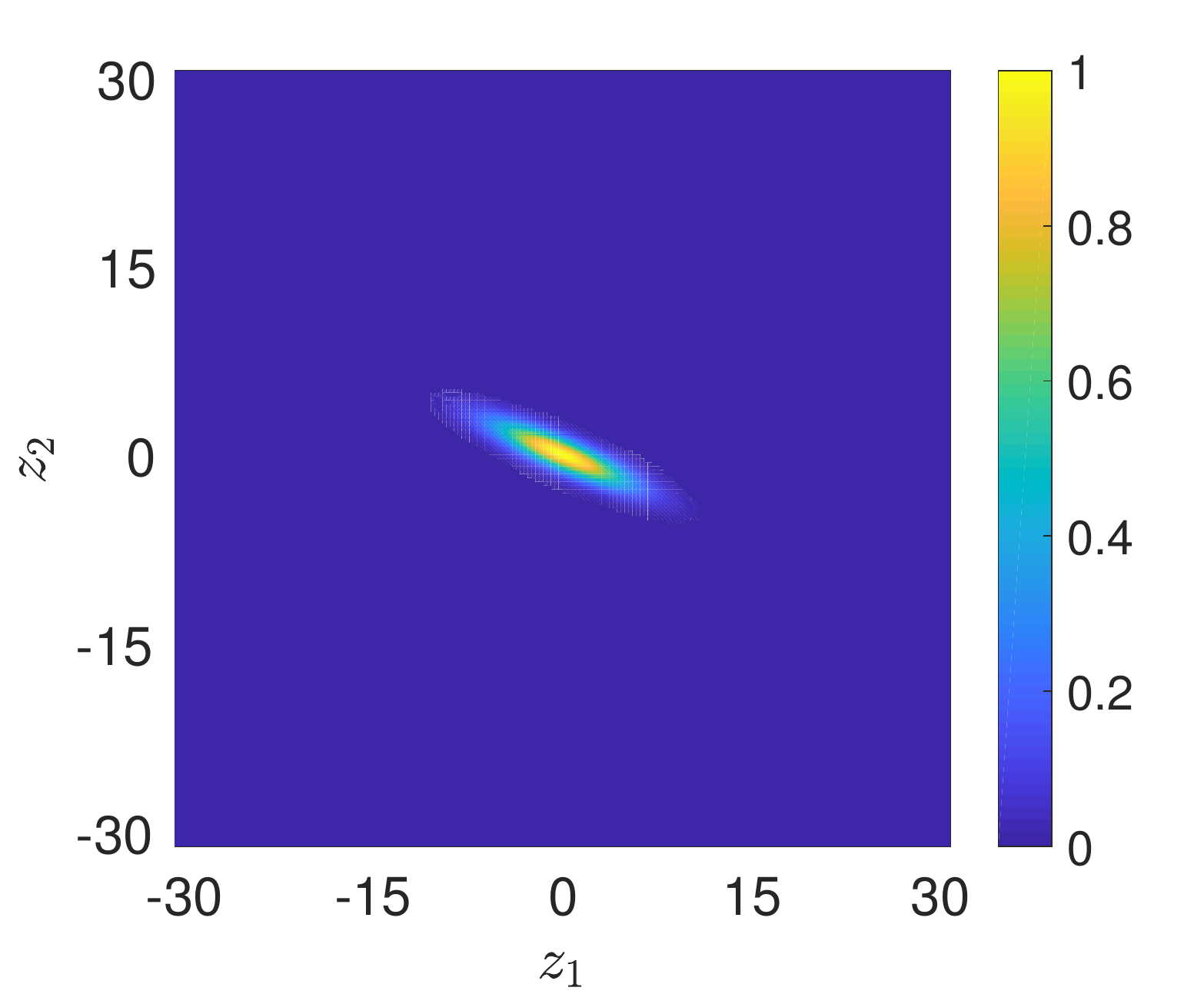}\hspace{-0.0cm}
\includegraphics[height=4.2cm]{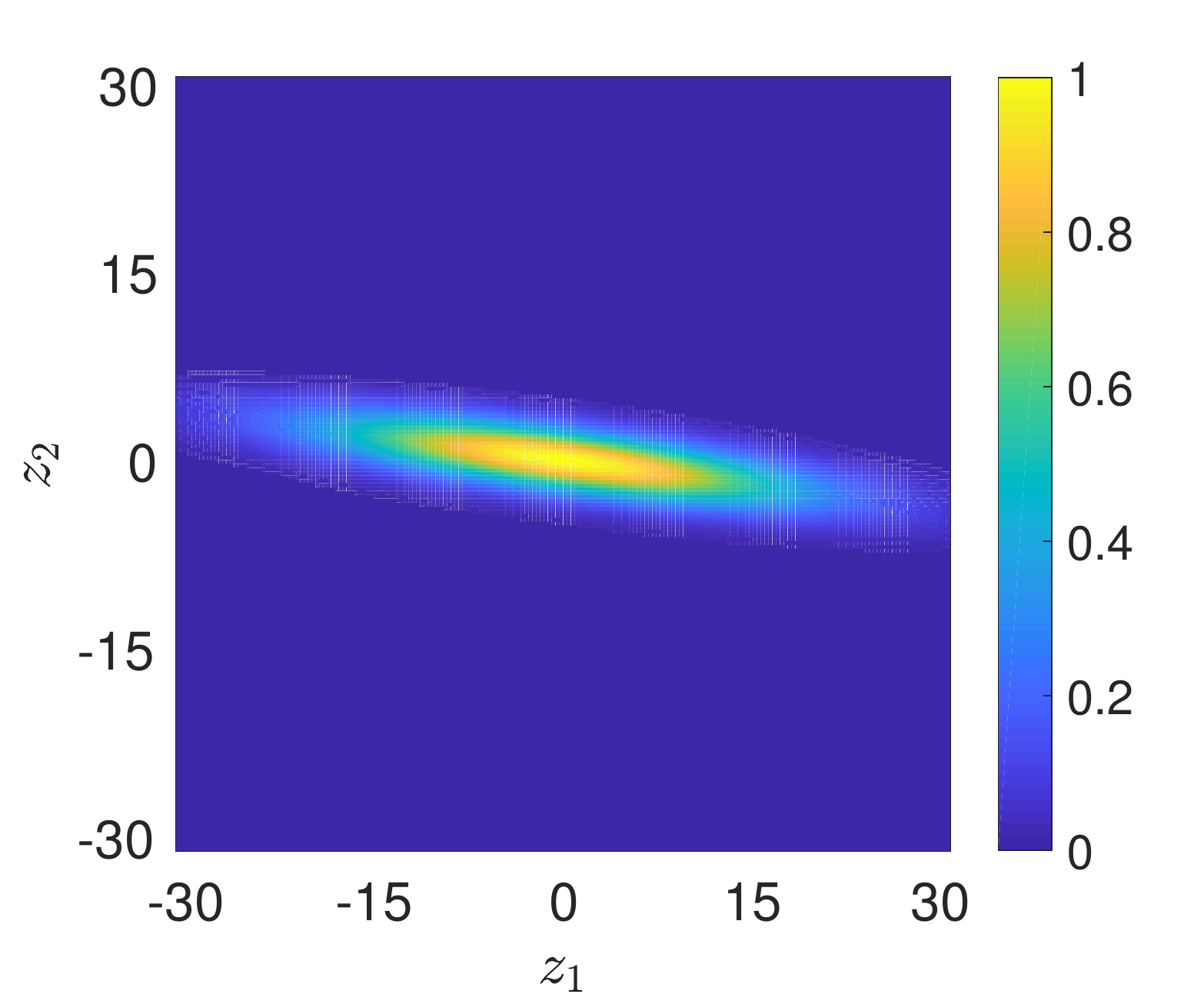}\hspace{-0.0cm}
}
\caption{
Solution to the multivariate PDE \eqref{PDE-advR} 
in $N=2$ dimensions. The stable spiral at the origin of the characteristic system attracts every curve in the phase space and ultimately yields $f=1$ 
everywhere after a transient.}
\label{fig:sol2D}
\end{figure}

Consider an initial-value problem
\begin{equation}
\frac{\partial f}{\partial t}+\sum_{j=1}^N\left(\sum_{k=1}^N 
C_{jk}z_k\right)\frac{\partial f}{\partial z_j}=0,
\qquad f(z_1,\ldots,z_N,0) = f_0(z_1,\ldots,z_N),
\label{PDE-advR}
\end{equation}
where $N$ is the number of independent variables forming the coordinate vector $\bm z = (z_1, \ldots, z_N)^\top$, and $C_{jk}$ are components of a given matrix of coefficients $\bm C$. 
The analytical solution to this initial-value problem is computed with the method of 
characteristics as~\cite{rhee2001}
\begin{align}
f(\bm z,t)=f_0 \left( \text{e}^{-t\bm C }\bm z\right).
\label{solutionF}
\end{align}
We choose the initial condition $f_0$ to be a fully separable product of exponentials,
\begin{align}
f_0(\bm z)=\text e^{-\left\|\bm z \right\|^2_2}.
\label{initialcondition}
\end{align}
The time dynamics of the solution \eqref{solutionF}-\eqref{initialcondition}
is entirely determined by the matrix $\bm C$.
In particular, if $\bm C$ has eigenvalues with 
positive real part then the matrix exponential 
$\exp(-t\bm C)$ is a contraction map. 
In this case, $f(\bm z,t)\rightarrow 1$ as $t\rightarrow \infty$,
at each point $\bm z\in \mathbb{R}^N$. 
Figure \ref{fig:sol2D} exhibits the analytical solution~\eqref{solutionF} 
generated by a $2\times 2$ matrix $\bm C$, whose eigenvalues are complex conjugate with positive real part.
The side length of the hypercube cube that encloses 
any level set of the solution at time $t$ depends on 
the number of dimensions $N$. In particular, the side-length 
of the hypercube that encloses 
the set $\{\bm z\in \mathbb{R}^N\,|\, f(\bm z,t)\geq \epsilon\}$
is given by 
\begin{equation}
\frac{1}{2}\sqrt{\frac{-\log(\epsilon)}{ N \lambda_\text{min}}}\leq b 
\leq \frac{1}{2}\sqrt{\frac{-\log(\epsilon)}{\lambda_\text{min}}}.
\label{est:b}
\end{equation}
where $\lambda_\text{min}$ is the smallest eigenvalue of the matrix 
$\exp( t \bm C)^\top \exp(t \bm C)$.
Equation \eqref{PDE-advR} can be written in the form 
\eqref{PDE0}, with separable operator $L$ (separation rank $N^2$). 
Specifically, the one-dimensional operators $L^l_j(z_j)$ appearing 
in~\eqref{Lz} are explicitly defined in Table~\ref{tab:ordering0}.

\begin{table}[htbp]
\centering
\begin{tabular}{c|ccccc}
$q$ & $\alpha_q$ & $L_1^q$ & $L_2^q$ & $\cdots$ & $L_N^q$\\
\hline\\
$1$ &  $-C_{11}$ &  $z_1\partial_{z_1}$ & $1$ & $\cdots$& $1$\\
$2$ &$-C_{12}$ &  $      \partial_{z_1}$ & $z_2$ & $\cdots$& $1$\\
$\vdots $&  $\vdots $ & $\vdots$ & $\vdots$ &$\vdots$ &$\vdots$ \\
$N$ &$-C_{1N}$ & $\partial_{z_1}$       & $1$ & $\cdots$& $z_N$\\
$N+1$ &$-C_{21}$ & $z_1$       & $\partial_{z_2}$ & $\cdots$& $1$\\
$N+2$ &$-C_{22}$ & $1$       & $z_2\partial_{z_2}$ & $\cdots$& $1$\\
$\vdots $&  $\vdots $ & $\vdots$ & $\vdots$ &$\vdots$ &$\vdots$ \\
$2N$ &$-C_{2N}$ & $1$       & $\partial_{z_2}$ & $\cdots$& $z_N$\\
$\vdots $&  $\vdots $ & $\vdots$ & $\vdots$ &$\vdots$ &$\vdots$ \\
$N^2-N+1$ &$-C_{N1}$ & $z_1$ & $1$ & $\cdots$& $\partial_{z_N}$\\
$N^2-N+2$ &$-C_{N2}$ & $1$ & $z_2$ & $\cdots$& $\partial_{z_N}$\\
$\vdots $&  $\vdots $ & $\vdots$ & $\vdots$ &$\vdots$ &$\vdots$ \\
$N^2$ &$-C_{NN}$ & $1$ & $1$ & $\cdots$& $z_N\partial_{z_N}$\\
\end{tabular}
\caption{Advection equation \eqref{PDE-advR}. Ordering 
of the linear operators defined in \eqref{Lz}.}
\label{tab:ordering0}
\end{table}

We used tensor methods, both canonical tensor decomposition and 
hierarchical tensor methods, and explicit time stepping (Section~\ref{sec:BoltzmanExplicit}) to solve numerically~\eqref{PDE-advR} in the 
periodic hypercube\footnote{Such domain is chosen large 
enough to accommodate periodic (zero) boundary 
conditions in the integration period of interest.} $[-60,60]^N$.  Each tensor component was discretized by using a Fourier spectral expansion with $600$ nodes in each variable (e.g., basis functions in \eqref{gfun} or \eqref{HTD}). 
The accuracy of the numerical solution was quantified in terms of the time-dependent relative error
\begin{equation}
\epsilon_m(t) = \left|\frac{f(\bm z^*,t)-\hat{f}(\bm z^*,t)}{f(\bm z^*,t)}\right|,
\label{relative_error}
\end{equation}  
where $f$ is the analytical solution \eqref{solutionF}, $\hat{f}$ 
is the numerical solution obtained by using the 
canonical or hierarchical tensor methods with separation rank $r$. 

\begin{figure}[htbp]
\centerline{\hspace{0.8cm}Canonical Tensor Method\hspace{3.0cm} Hierarchical Tensor Method}
\noindent
\centerline{
\rotatebox{90}{\hspace{2cm}two dimensions}\hspace{0.2cm}
\includegraphics[height=5.5cm]{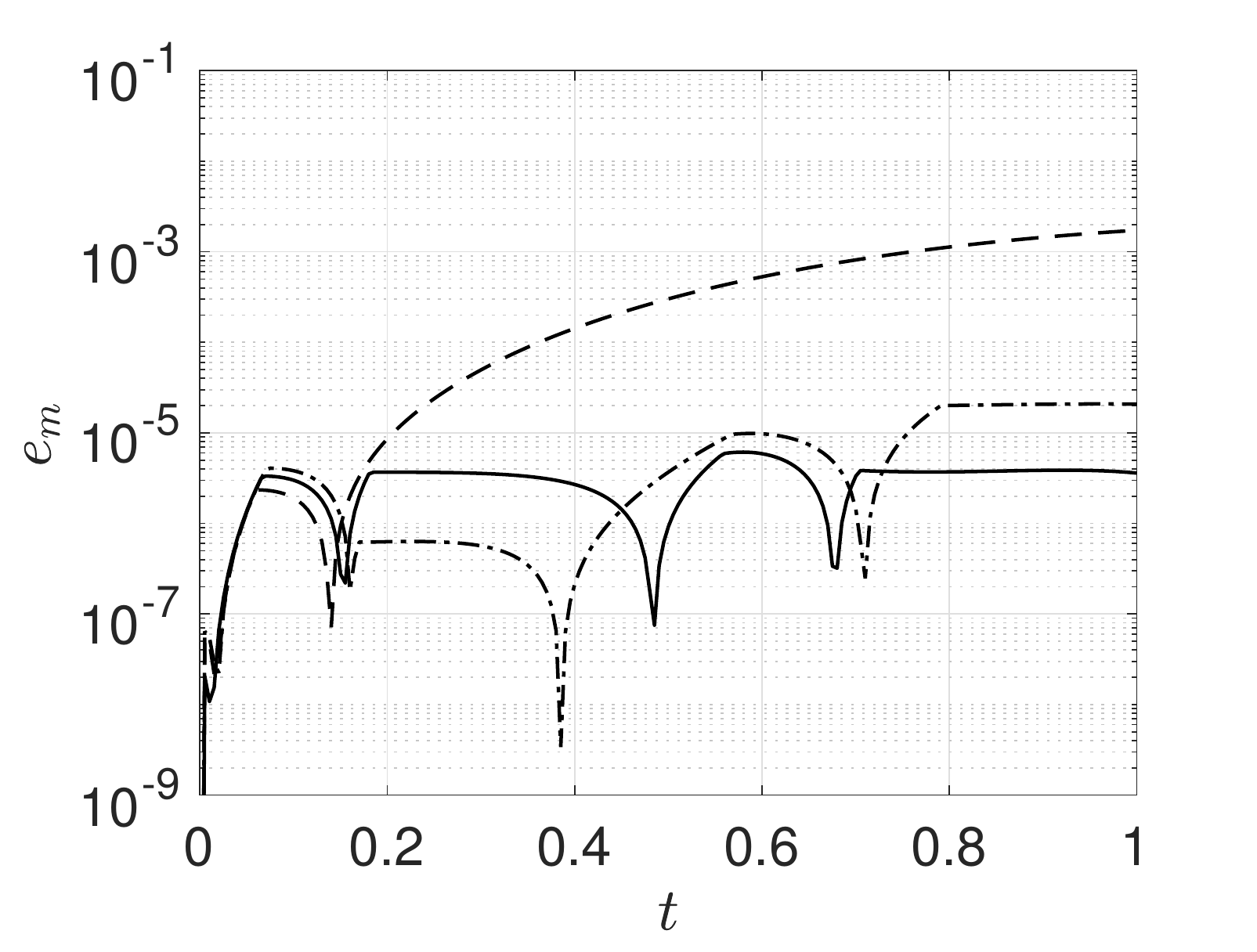}
\includegraphics[height=5.5cm]{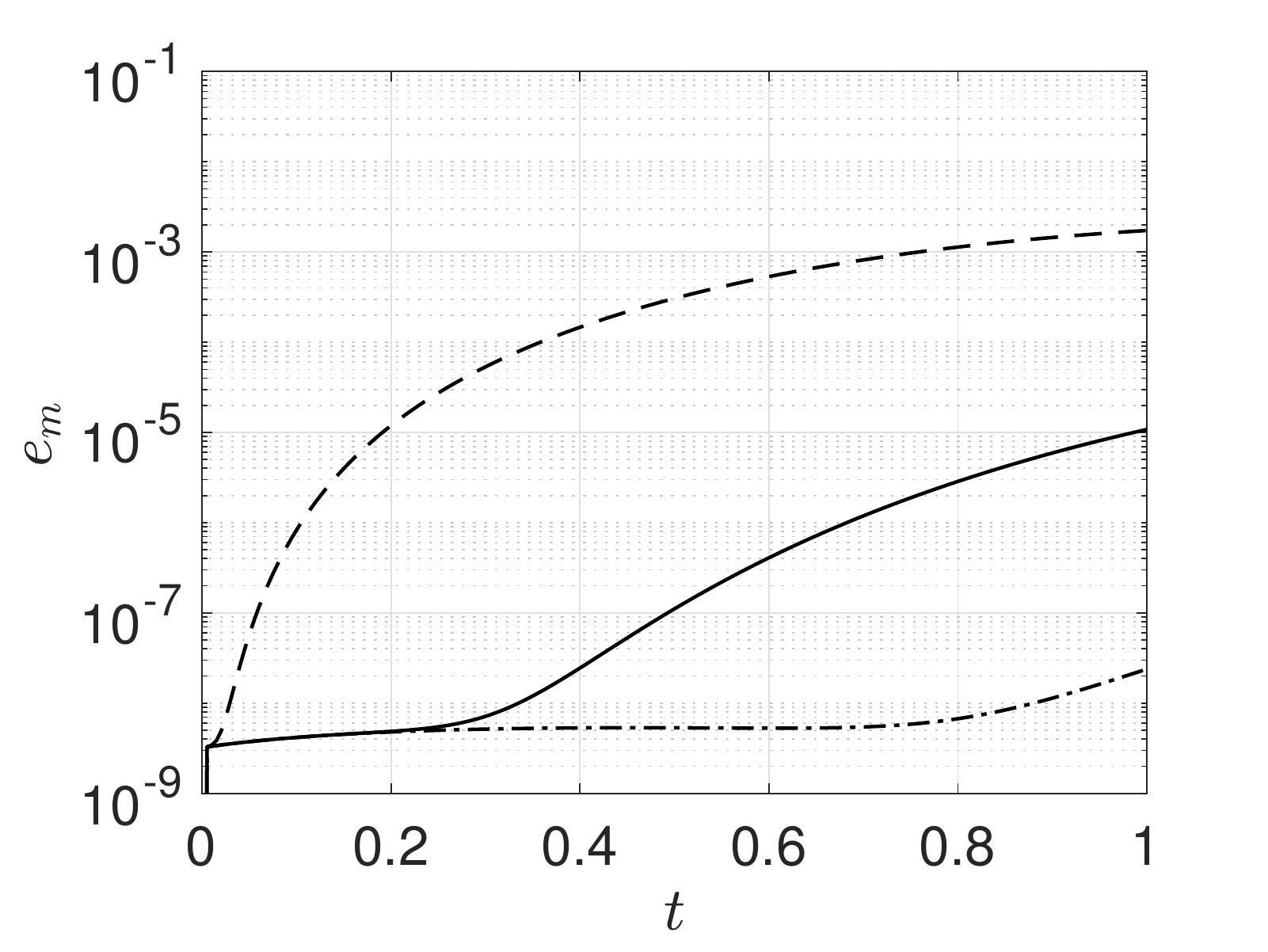}
}
\noindent
\centerline{
\rotatebox{90}{\hspace{2cm}three dimensions}\hspace{0.2cm}
\includegraphics[height=5.5cm]{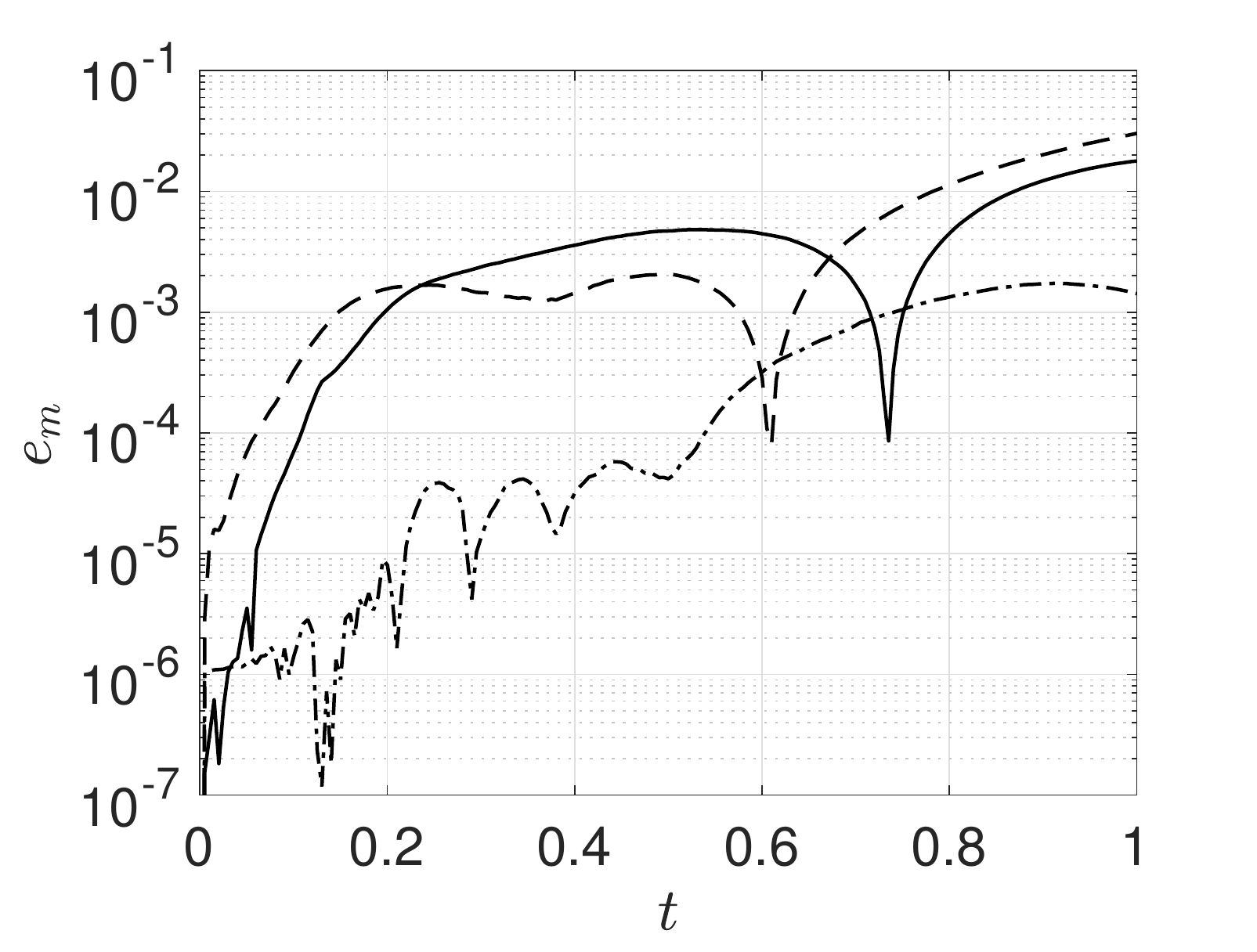}
\includegraphics[height=5.5cm]{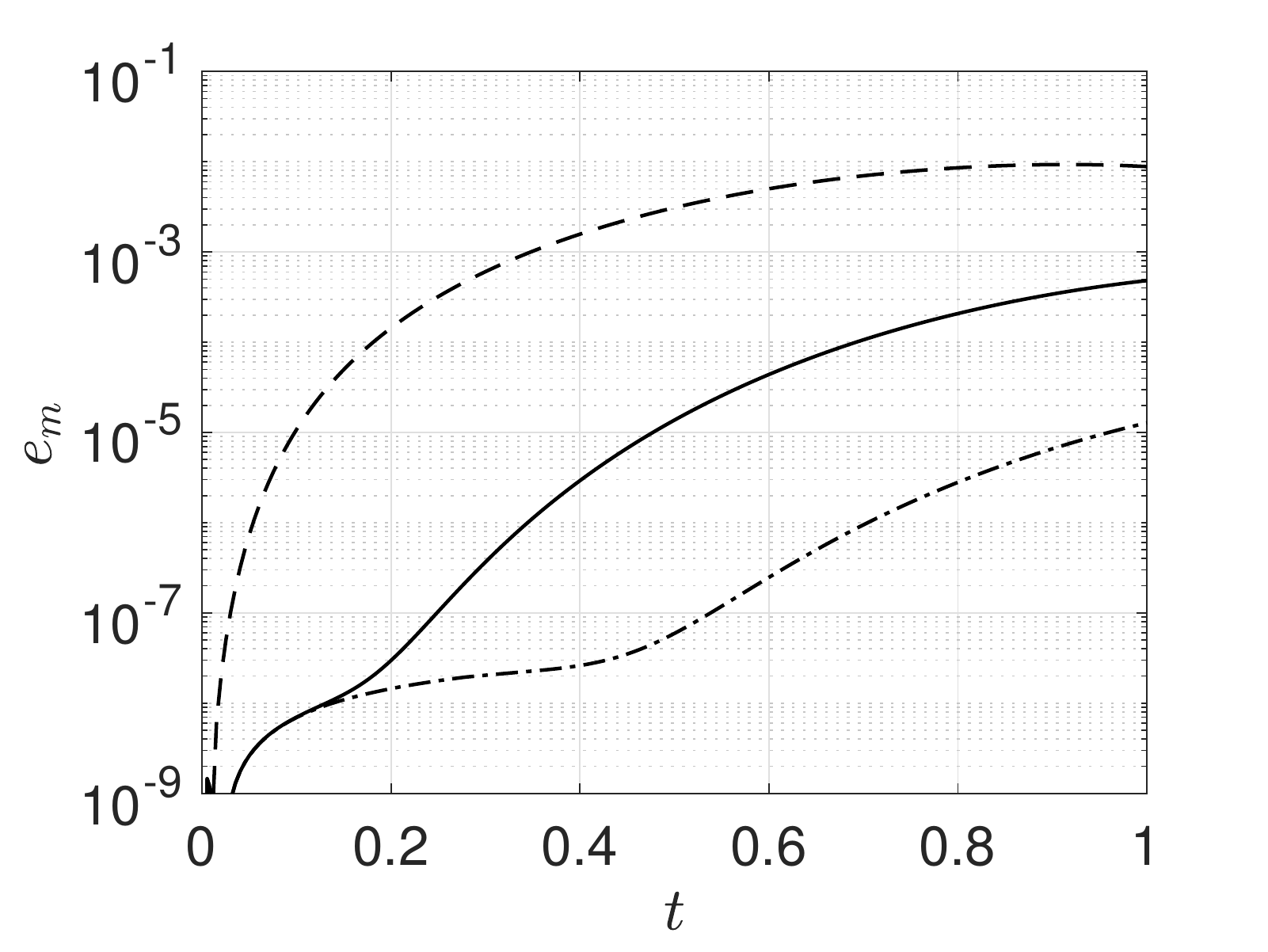}
}

\centerline{
\rotatebox{90}{\hspace{2cm}six dimensions}\hspace{0.2cm}
\includegraphics[height=5.5cm]{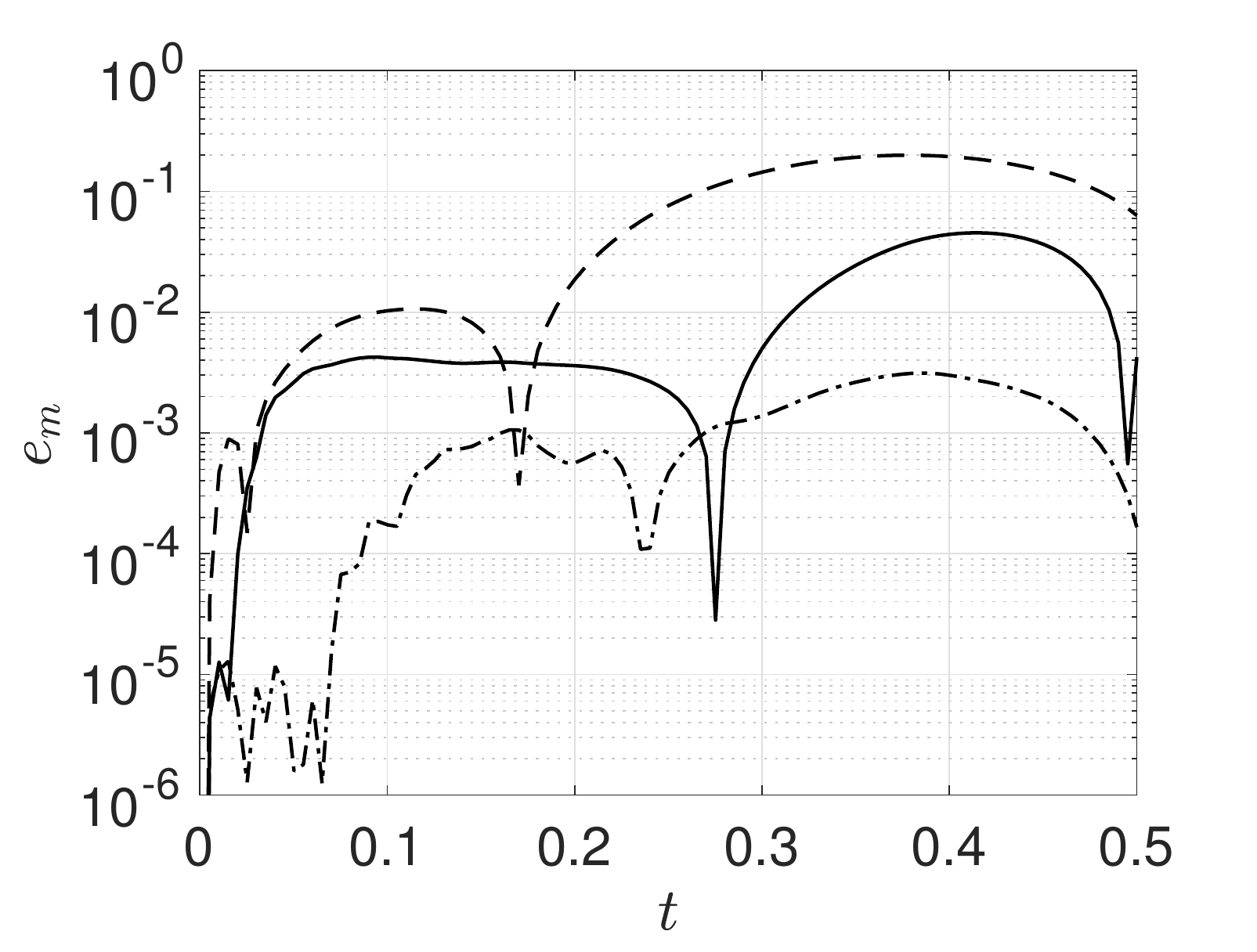}
\includegraphics[height=5.5cm]{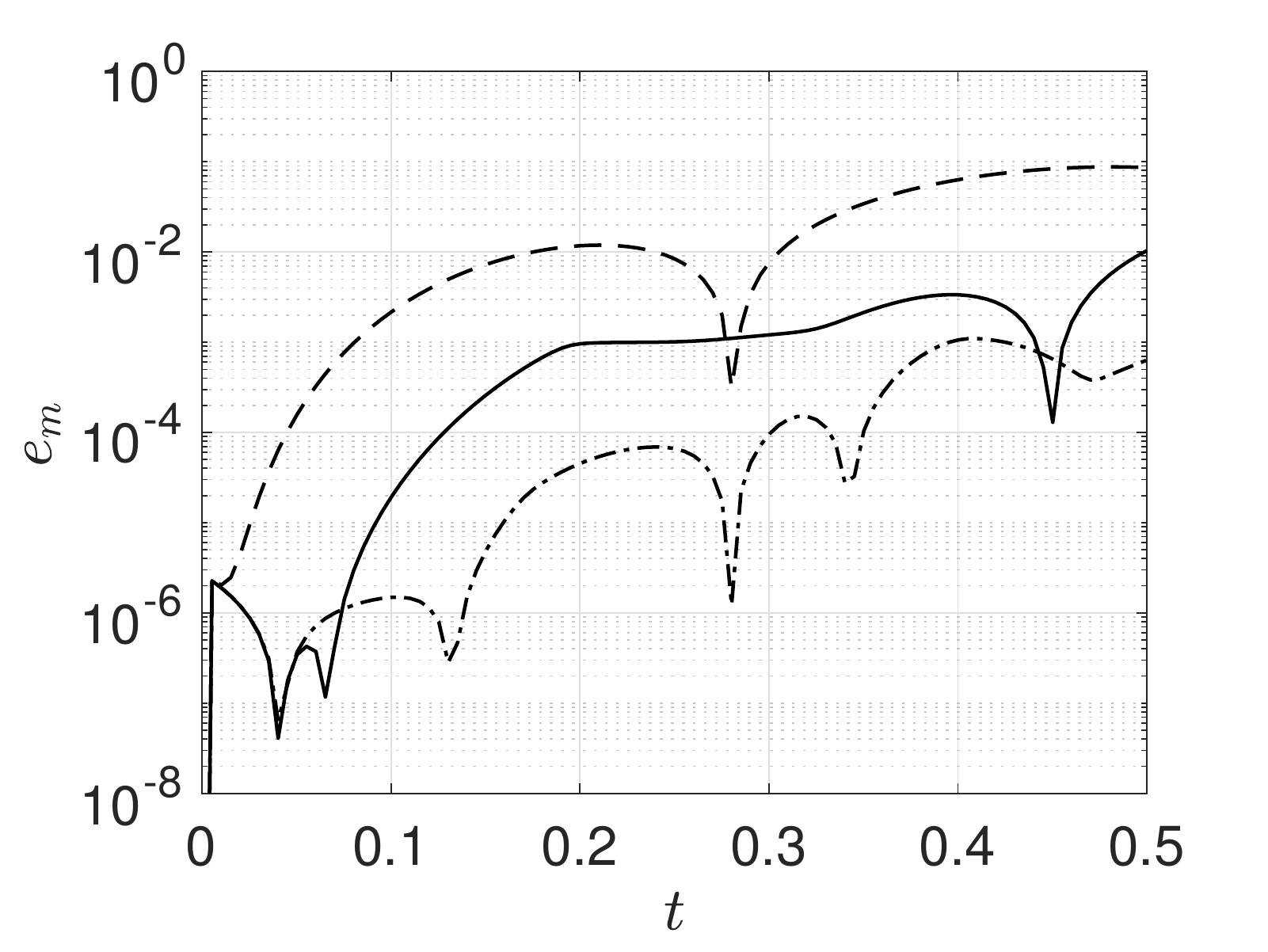}
}
\centerline{$--$ $r_{max}=4$\hspace{1cm} --- $r_{max}=8$
\hspace{1cm} $-\cdot -$ $r_{max}=12$
}
\caption{Advection equation \eqref{PDE-advR}. 
Relative pointwise error \eqref{relative_error} for 
several values of the separation ranks $r$ and the number 
of independent variables $N$. These results are obtained 
by using the canonical and hierarchical tensor methods with 
explicit time stepping (see Section \ref{sec:BoltzmanExplicit}).
The separation rank is computed adaptively 
at each time step up to the maximum value $r_\text{max}$.}
\label{fig:ALS-CP-HT}
\end{figure}

Figure \ref{fig:ALS-CP-HT} shows the relative pointwise error \eqref{relative_error}, computed at $\bm z^*=(h,h,\ldots,h)$ with $h= 0.698835274542439$, for different separation ranks $r$ and  
different number of independent variables $N$. 
As expected, the accuracy of the numerical solution increases with the separation rank $r$. 
Also, the relative error increases with the number of dimensions $N$.
It is worthwhile emphasizing that the rank reduction in 
canonical tensor decomposition at each time step is based 
on a randomized algorithm that requires 
initialization at each time-step. This means that 
results of simulations with the same nominal 
parameters may be different. On the other hand, tensor
methods based on a hierarchical singular 
value decomposition, such as the hierarchical Tucker 
decomposition \cite{grasedyck2010,grasedyck2017}, do 
not suffer from this issue.
The separation rank of both canonical and hierarchical tensor 
methods is computed adaptively up 
to the maximum value $r_\text{max}$ specified 
in the legend of Figure~\ref{fig:ALS-CP-HT}. 
In two dimensions, the CP-ALS algorithm 
yields the similar error plots for 
$r_\text{max}=8$ and $r_\text{max}=12$. This is because $r < 8$ in both cases, and 
throughout the simulation up to $t=1$. 
The difference between these results is due to the 
random initialization required by the ALS algorithm 
at each time step.

\subsection{Boltzmann-BGK Equation}
\label{sec:BoltzmanDiscretization}
In this section we develop an accurate ALS-based algorithm  to solve  a linearized BGK equation (see Appendix~\ref{sec:Boltzmann} for details of its derivation),
\begin{equation}
\frac{\partial f}{\partial t}+\bm v \cdot \nabla_{\bm x} f = 
\nu I, \qquad I = f_\text{eq}(\bm v)-f. 
\label{BGK_LIN}
\end{equation}
Here, $f(\bm z,t)$ is a probability density function is six phase variables plus 
time; the coordinate vector $\bm z = (\bm x, \bm v)^\top$ is composed of three spatial dimensions $\bm x = (x_1,x_2,x_3)^\top$ and three components of the velocity vector $\bm v = (v_1,v_2,v_3)^\top$; $f_\text{eq}$ denotes a locally Maxwellian distribution,
\begin{equation}
f_\text{eq}(\bm v) = \frac{\rho}{(2\pi R T)^{3/2}}
\exp\left(-\frac{\left\|\bm v\right\|_2^2}{2R T}\right),
\label{EqMaxwell}
\end{equation}
with uniform gas density $\rho$, temperature $T$, and velocity $\bm v$; $R$ is the gas constant; and $\nu = \kappa \rho T^{1-\mu}$ with positive constants $\kappa$ and $\mu$. The solution to the linearized Boltzmann-BGK equation \eqref{PDE11} is also
computable analytically, which provides us with a benchmark solution 
to check the accuracy of our algorithms. 

As before we start by rewriting this equation 
in an operator form, 
\begin{equation}
\frac{\partial f}{\partial t} = L f + C
\label{PDE11} 
\end{equation}
where 
\begin{align}\label{opL}
L(\bm x,\bm v) = -\bm v\cdot \nabla_{\bm x}- \nu I,\qquad
C(\bm v) = \nu f_\text{eq}(\bm v).
\end{align}
Note that $L$ is a {\em separable} linear 
operator with separation rank $r_L=4$, which can be rewritten in the form 
\begin{equation}
L(\bm z)= \sum_{q=1}^{4} \alpha_q L^q_1(z_1)\cdots L_6^q(z_6),
\label{eq:operators}
\end{equation}
for suitable one-dimensional linear operators $L^q_j(a_j)$ 
defined in Table \ref{tab:ordering}. Similarly, $C(\bm v)$ in \eqref{opL} is separated as 
\begin{equation}
C(\bm v) = C_1(v_1)C_2(v_2)C_3(v_3), \quad \textrm{where}\quad 
C_i(v_i)=\frac{(\nu \rho)^{1/3}}{(2\pi R T)^{1/2}}
\exp\left(-\frac{v_i^2}{2R T}\right).
\label{sepC}
\end{equation}

\begin{table}[htbp]
\centering
\begin{tabular}{c|cccccccc}
$q$ & $\alpha_q$ & $L_1^q$ & $L_2^q$ & $L_3^q$ & $L_4^q$ & $L_5^q$ & $L_6^q$\\
\hline
$1$ &  $-\nu$ &  $1$ & $1$ & $1$ & $1$& $1$ & $1$ \\
$2$ &  $-1$ &  $\partial/\partial x_1$ & $1$ & $1$ & $v_1$& $1$ & $1$ \\
$3$ &  $-1$ &  $1$ & $\partial/\partial x_2$ & $1$ & $1$& $v_2$ & $1$ \\
$4$ &  $-1$ &  $1$ & $1$ & $\partial/\partial x_3$ & $1$& $1$ & $v_3$ \\
\end{tabular}
\caption{Ordering of the linear operators defined 
in  equation~\eqref{eq:operators}.}
\label{tab:ordering}
\end{table}

From a numerical viewpoint, the solution to \eqref{PDE11} can be 
represented by using any of the tensor series expansion
we discussed in Section \ref{sec:tensorseries}. In particular, 
hereafter we develop an algorithm based on canonical tensor tensor 
decomposition, alternating least squares, and implicit time 
stepping (see Section \eqref{sec:BoltzmanImplicit}).
To this end, we discretize \eqref{PDE11} in time with 
the Crank-Nicolson method to obtain
\begin{equation}
\left[I-\frac{\Delta t}{2} L \right]f_{n+1}=
\left[I+\frac{\Delta t}{2} L \right]f_n+\Delta t C + \Delta t\tau_n, 
\label{CNFDE1}
\end{equation}
where $f_n(\bm z)= f(\bm z,t_n)$ and $\tau_n$ it the local truncation 
error of the Crank-Nicholson method at $t=t_n$ (\cite{quarteroni2007}, p. 499). 
This equation can be written in a compact notation as 
\begin{equation}
Af_{n+1}=Bf_n + \Delta t C +\Delta t \tau_n,
\end{equation}
where  $A$ and $B$ are separable operators in the form \eqref{BLs}, 
where all quantities are defined in Table \ref{tab:AL}.

\begin{table}[htbp]
\centering
\begin{tabular}{c|ccccccccc}
$q$ &  $\eta_q$ & $\zeta_q$ & $E_{1}^q$ & $E_{2}^q$ & $E_{3}^q$ & $E_{4}^q$ & $E_{5}^q$ & $E_{6}^q$\\
\hline\\
$1$ &  $\nu (1+\Delta t/2)$ &  $\nu (1-\Delta t/2)$ &  
$1$ & $1$ & $1$ & $1$& $1$ & $1$\\
$1$ &  $\Delta t/2$ &  $-\Delta t/2$ &
$\partial/\partial x_1$ & $1$ & $1$ & $v_1$& $1$ & $1$\\
$1$ &  $\Delta t/2$ &  $-\Delta t/2$ &
$1$ & $\partial/\partial x_2$ & $1$ & $1$& $v_2$ & $1$\\
$1$ &  $\Delta t/2$ &  $-\Delta t/2$ &
$1$ & $1$ & $\partial/\partial x_3$ & $1$& $1$ & $v_3$
\end{tabular}
\caption{Ordering of the linear operators $A$ and $B$ defined 
in  \eqref{BLs}.}
\label{tab:AL}
\end{table}
A substitution of the canonical tensor  
decomposition\footnote{Recall that the 
functions $G_{k}^l(z_k,t_n)$ 
are in the form 
\begin{equation}
G_{k}^l(z_k,t_n)=\sum_{s=1}^Q 
\beta^l_{ks}(t_n)\phi_s(z_k),
\end{equation}
$\beta^l_{ks}(t_n)$ ($l=1,...,r$, $k=1,...,6$, $s=1,...,Q$) 
being the degrees of freedom. }
\begin{equation}
\hat{f}_{n+1} = \sum_{l=1}^r \prod_{k=1}^6 G^l_k(z_k,t_{n+1})
\label{CP-numerical}
\end{equation}
into equation \eqref{CNFDE} yields the residual 
\begin{equation}
R(\bm z,t_{n+1}) = A(\bm z)\hat{f}_{n+1}(\bm z)-B(\bm z)\hat{f}_n(\bm z) + \Delta t C(z_4,z_5,z_6),
\label{rRr}
\end{equation}
which can be minimized by using the alternating least squares method, as
we described in Section \ref{sec:BoltzmanImplicit} 
to obtain the solution at time $t_{n+1}$.

\paragraph{Nonlinear Boltzmann-BGK Model}

The numerical tensor methods we discussed in Section \ref{sec:BoltzmanImplicit} and Section \ref{sec:BoltzmanExplicit} can be extended 
to compute the numerical solution of the fully nonlinear Boltzmann-BGK model
\begin{equation}
\frac{\partial f}{\partial t} +\bm v\cdot \nabla_{\bm x} f = 
 \nu(\bm x,t)\left(f_{eq}(\bm x,\bm v,t)-f(\bm x,\bm v,t)\right).
 \label{fullBGK}
\end{equation}
Here, $f_{eq}(\bm x,\bm v,t)$ is the equilibrium distribution \eqref{LocalMaxwell}, 
while the collision frequency $\nu(\bm x,t)$ can can be expressed, e.g., 
as in  \eqref{collfreq}. From a mathematical viewpoint both quantities 
$f_{eq}$ and $\nu$ are nonlinear functionals of the PDF $f(\bm x,\bm v,t)$ 
(see equations \eqref{rho}-\eqref{T}). Therefore \eqref{fullBGK} is  
an advection equation driven driven by a nonlinear functional 
of $f(\bm x,\bm v,t)$. The evaluation of the integrals appearing 
in \eqref{rho}-\eqref{T} poses no great challenge if we represent $f$ as  
a canonical tensor \eqref{functional-SSE} or as a hierarchical Tucker tensor
\eqref{HTD}. Indeed, such integrals can be factored as sums of 
product of one-dimensional integrals in a tensor representation. Moreover, 
computing the inverse of $\rho(\bm x,t)$ in 
\eqref{u} and \eqref{T} is relatively simple in a collocation 
setting. In fact, we just need to evaluate the full tensor $\rho(\bm x,t)$, on 
a three-dimensional tensor product grid and 
compute  $1/\rho(\bm x,t)$ pointwise. If needed, we could then compute 
the tensor decomposition of $1/\rho(x,t)$, e.g., by using the canonical tensor representations and the ALS algorithm we discussed in Section \ref{sec:CP}.
This allows us to evaluate the BGK collision operator  and represent 
it in a tensor series expansion. The solution to the Boltzmann BGK 
equation can be then computed by operator splitting methods \cite{dimarco2018,desvillettes1996}.

\paragraph{Initial and Boundary conditions}
To validate the proposed alternating least squares algorithm we 
set the initial condition to be either coincident with the homogeneous
Maxwell-Boltzmann distribution \eqref{EqMaxwell}, i.e.,
\begin{equation}
  f \left(\bm x,\bm v,0\right) =
  \frac{\rho}{\left(2 \pi R T\right)^{3/2}} 
  \exp\left(-\frac{\left\|\bm v\right\|_2^2}{2 R T}\right)
  \label{IC}
\end{equation}
or a slight perturbation of it. 
The PDF \eqref{IC} is obviously separable with separation rank one. 
Moreover, we set the computational domain to be a
the periodic hyperrectangle $[-b_x,b_x]^3\times[-b_v,b_v]^3$.
It is essential that such domain is chosen large
enough to accommodate the support of the PDF at all times, thus
preventing physically unrealistic correlations between the distribution function. 
More generally, it is possible to enforce other types of boundary conditions 
by using appropriate trial functions, or by a mixed approach
\citep{shuleshko1959,snyder1964,zinn1968} in the case of 
Maxwell boundaries.

Unless otherwise stated, all simulation parameters  
are set as in Table \ref{tab:param}. Such setting corresponds 
to the problem of computing the dynamics of Argon within 
the periodic hyperrectangle $[-b_x,b_x]^3\times [-b_v,b_v]^3$. 
We are interested in testing our algorithms in two different regimes: 
i) steady state, and ii) relaxation to statistical equilibrium. 
\begin{table}[t]
\centering 
\begin{tabular}{|l|l|l|}
\hline
Parameter             & Symbol                 & value                               \\                                                                           
\hline                                                                               
Temperature           & $T$                    & $300\, \si{\K}$                       \\
Number density        & $n$                    & $2.4143 \cdot 10^{25}\, \si{\m^{-3}}$ \\
Specific gas constant & $R$                & $208 \,\si{\J.\kg^{-1}.\K^{-1}}$      \\
Relaxation time       & $\tau_{R}$             & $0.40034\, \si{\s}$                   \\
Collision frequency   & $\nu$                  & $1/\tau_{R} $           \\
Position domain size  & $b_{x}$                & $500\, \si{\m}$                      \\
Velocity domain size  & $b_{v}$                & $5 \sqrt{R_{S} T}$ \\
Time step             & $\Delta t$             & $0.01 \tau_{R} $             \\
Number of iterations  & $N_{\txt{Iter}}$       & $1000$                              \\
ALS tolerance & $\epsilon_{\txt{Tol}}$ & $10^{-8}$                           \\
Series truncation     & $Q$                    & $11$                                \\ 
\hline
\end{tabular}
\center
\caption{Simulation parameters we employed in the numerical approximation 
of the Boltzmann-BGK model. These parameters correspond to a 
kinetic simulation of dynamics of Argon in the periodic hyperrectangle 
$[-b_x,b_x]^3\times [-b_v,b_v]^3$. The time step in the implicit 
Crank-Nicolson method is $\Delta t$, while the tolerance 
in the ALS iterations at each time step is $\epsilon_{\txt{Tol}}$. 
}
\label{tab:param}
\end{table}

\subsubsection{Canonical Tensor Decomposition of the Maxwell-Boltzmann Distribution}
\begin{figure}[t]
\center
\subfloat{
\includegraphics[width=0.50\textwidth]{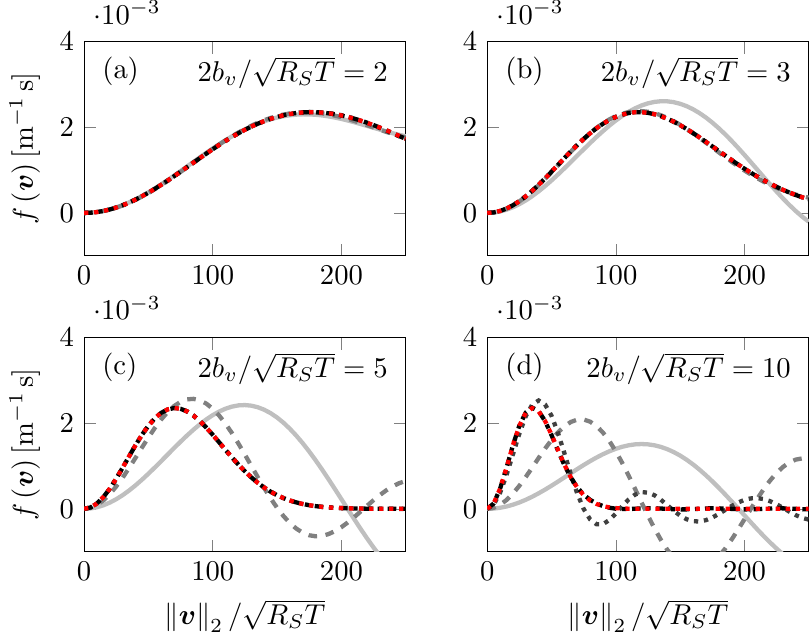}
}
\subfloat{
\includegraphics[width=0.50\textwidth]{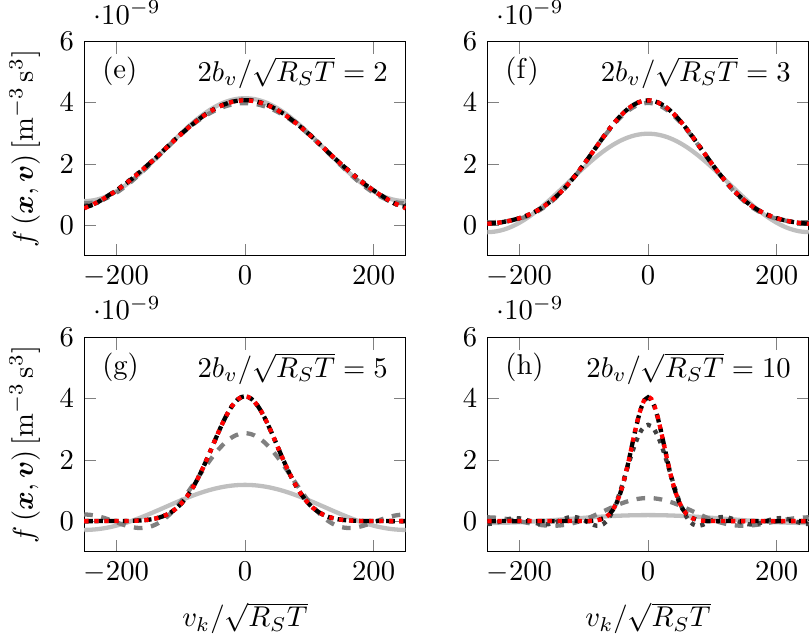}
} \\
\caption{One-particle Maxwell-Boltzmann distribution $f_{eq}(\bm v)$ 
(\protect\includegraphics{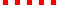}) as a function of the 
dimensionless molecular
speed $\left\|\bm v\right\|_2/\sqrt{R T}$ for different 
dimensionless domain sizes $b_{v}/\sqrt{R T}$.  
In each case, we demonstrate convergence of the canonical tensor 
decomposition \eqref{functional-SSE} of $f_{eq}$ 
as we increase the number of basis functions
$Q$.  Specifically, we plot the cases: $Q=3$ (\protect\includegraphics{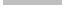}), 
$Q=5$ (\protect\includegraphics{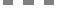}),  $Q=11$ (\protect\includegraphics{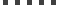}) and $Q=19$ (\protect\includegraphics{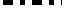}).
}

\label{fig:pdfQ}
\end{figure}

We begin with a preliminary convergence study of the canonical tensor 
decomposition \eqref{functional-SSE} applied to the 
one-particle Maxwell distribution $f_{eq}(\bm v)$. Such study 
allows us to identify the size of the hyperrectangle that 
includes the support of the Maxwellian equilibrium  distribution. 
This is done in Figure \ref{fig:pdfQ}, where we plot 
the tensor expansion of $f_{eq}(\bm v)$ as function 
of the non-dimensional molecular speed $\left\|\bm v\right\|_2/\sqrt{R_{S} T}$, 
for different number of basis functions $Q$ in \eqref{gfun}, and for
various values of the dimensionless velocity domain size $b_{v}/\sqrt{R T}$. 
It is seen that for small values of $b_{v}/\sqrt{R T}$, a small number of 
basis functions $Q$ is sufficient to capture the exact equilibrium 
distribution. At the same however, small values of
$b_{v}/\sqrt{R T}$ result in distribution functions with 
truncated tails, as can be seen in Figures \ref{fig:pdfQ} (a), \ref{fig:pdfQ}
(b), and \ref{fig:pdfQ} (e). 
Figures~\ref{fig:pdfQ} (c) and \ref{fig:pdfQ} (g)
demonstrate that the dimensionless velocity domain size $b_{v}/\sqrt{R  T} = 5$ 
is sufficient to avoid truncation of the tail of the distribution, 
and that a series expansion with $Q=11$ in each variable is sufficient 
to accurately capture the equilibrium distribution. This justifies 
the choice of parameters in Table \eqref{tab:param}.

\begin{figure}[htpb]
\centering
\subfloat{
\includegraphics[width=0.50\textwidth]{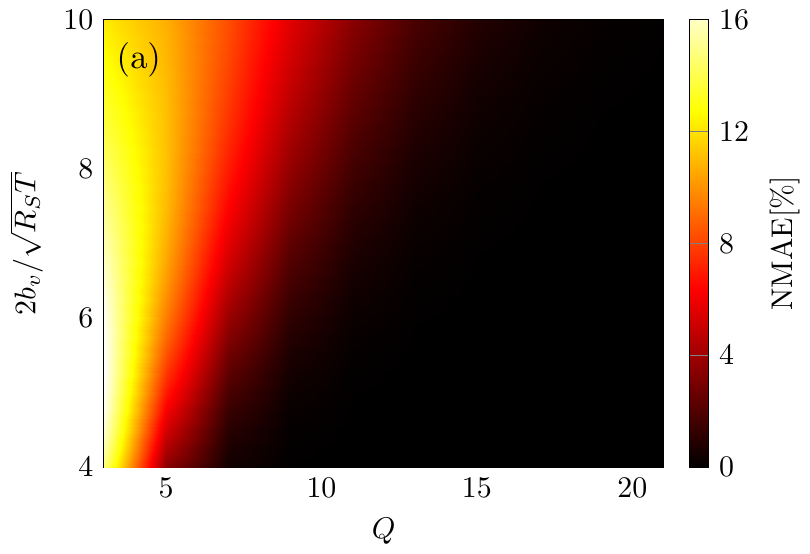}
}
\subfloat{
\includegraphics[width=0.50\textwidth]{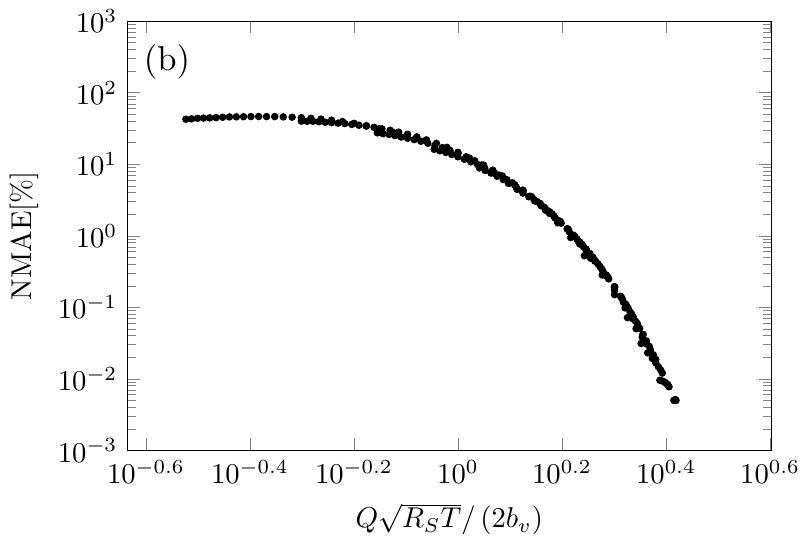}
} \\
\caption{(a) Normalized Mean Absolute Error (NMAE) \eqref{NMAE} 
between the analytical and the canonical tensor series 
expansion of the Maxwell-Boltzmann equilibrium distribution 
versus versus $Q$ and the dimensionless velocity domain 
size $b_{v}/\sqrt{R T}$. In Figure (b) we plot the NMAE
in a log-log scale as function of $Q \sqrt{R T}/b_{v}$. We emphasize that 
the range of $Q$ and $b_{v}/\sqrt{R  T}$ is the same in both Figures.}
\label{fig:NMAE}
\end{figure}

A more detailed analysis of the effects of the truncation order $Q$
and the dimensionless velocity domain size $b_{v}/\sqrt{R T}$ 
on the approximation error in $f_{eq}$ is done in Figure \ref{fig:NMAE}(a), 
where we plot the Normalized Mean Absolute Error (NMAE) 
versus $b_{v}/\sqrt{R T}$ and $Q$. The NMAE between two 
vectors $\bm X$ and $\bm Y$ of 
size $N$ is defined as: 
\begin{equation}
  \txt{NMAE} 
= \frac{1}{N}
  \frac{\left\|\bm X -\bm Y\right\|_1}{\max{(\bm X)}-\min{(\bm Y)}},
  \label{NMAE}
\end{equation}
where $\left\|\cdot \right\|_1$ is the standard vector 1-norm. 
In Figure  \ref{fig:NMAE}(b) we plot the NMAE as a function of 
$Q \sqrt{R T}/b_{v}$. This results in a collapse of the data which suggests 
that if one wants to double the non-dimensional velocity domain size while
maintaining the same accuracy, the series truncation order $Q$ needs 
to be doubled as well. In addition, it can be seen that 
starting at $Q \sqrt{R T}/b_{v} \approx 1$ the error decays 
rapidly with increasing $Q \sqrt{R T}/b_{v}$.

\subsubsection{Steady State Simulations}
\begin{figure}[t]
\center
\includegraphics[width=0.8\textwidth]{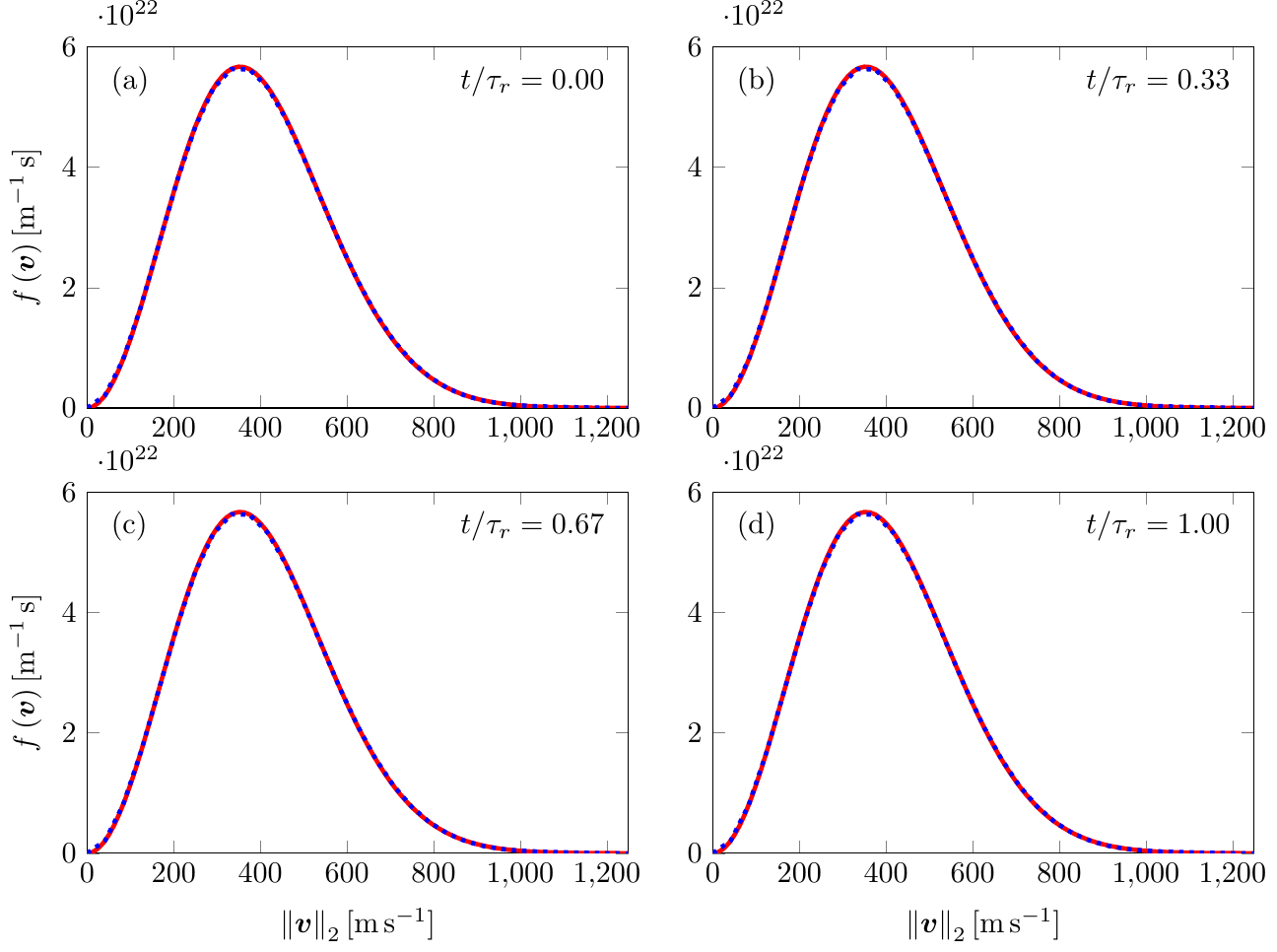}
\caption{Solution to the Boltzmann-BGK equation \eqref{BGK_LIN} 
at different dimensionless times $t/\tau_{R}$  
 ($\tau_R$ is 
the relaxation time in Table \ref{tab:param}). Specifically, 
we plot the Maxwell-Boltzmann equilibrium distribution 
(\protect\includegraphics{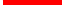}), its canonical tensor decomposition 
with $Q=11$ modes, and 
its the time evolution with the ALS algorithm we described in 
Section \ref{sec:BoltzmanImplicit}  (\protect\includegraphics{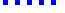}). 
It is seen that ALS has a stable fixed point at $f_{eq}(\bm v)$. This is 
an important test as, in general, ALS iterations are not granted to converge 
(e.g., \cite{uschmajew2012}).}
\label{fig:cConst}
\end{figure}

We first test the the Boltzmann-BGK solver we developed 
in Section \eqref{sec:BoltzmanImplicit} on an initial value problem 
where the initial condition is set to be the Maxwell-Boltzmann
distribution $f_{eq}(\bm v)$ (see equation \eqref{EqMaxwell}). 
A properly working algorithm should keep such equilibrium distribution 
unchanged as the simulation progresses.
In Figure~\ref{fig:cConst} we plot the 
distribution of molecular velocities $f(\bm v,t)$ as function of the molecular 
speed $\left\|\bm v\right\|_2\abs{\vec{\xi}}$ at various dimensionless 
times $t/\tau_R$, where $\tau_R$ is the relaxation time in Table \ref{tab:param}. 
It is seen that the alternating least squares algorithm we developed in 
Section \eqref{sec:BoltzmanImplicit} has a stable 
fixed point at the Maxwell-Boltzmann equilibrium 
distribution $f_{eq}$. This is an important test, as 
convergence of alternating least square is, in general, not granted 
for arbitrary residuals (see, e.g., \cite{uschmajew2012,comon2009}).

\begin{figure}[t]
\centering
\subfloat{
\includegraphics[width=0.30\textwidth]{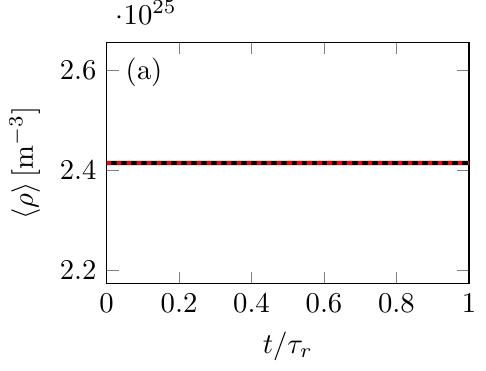}
}
\subfloat{
\includegraphics[width=0.30\textwidth]{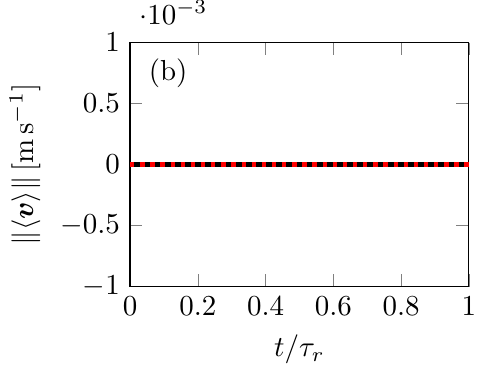}
}
\subfloat{
\includegraphics[width=0.30\textwidth]{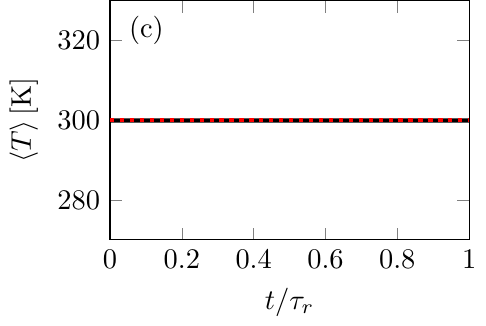}
} \\
\caption{Average density $\left<\rho\right>(t)$, 
absolute average velocity $\left\|\left< u\right>\right\|_2$, and 
average temperature  $\left<T\right> $ as function of the dimensionless
time $t/\tau_{R}$: simulation results (\protect\includegraphics{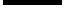}); 
results at equilibrium  (\protect\includegraphics{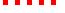}). 
The simulation results are obtained by setting the initial condition 
 in \eqref{BGK_LIN} equal to the Maxwell-Boltzmann equilibrium distribution 
 $f_{eq}$ and solving the Boltzmann-BGK model with the ALS algorithm described in Section \ref{sec:BoltzmanImplicit}. 
Clearly, ALS iterations  preserve the 
average density, velocity and temperature.}
\label{fig:nUTConst}
\end{figure}

We also computed the moments of the distribution 
function $f(\bm x,\bm v,t)$ at different times to verify 
whether ALS iterations preserve  (density), 
momentum (velocity), and kinetic energy (temperature). 
In particular, we computed  
\begin{align}
\left< \rho\right>(t)
&= 
  \frac{1}{b_{x}^{3}}
  \int_{-b_{x}}^{b_{x}}\cdots
  \int_{-b_{v}}^{b_{v}}
  f\left(\bm x,\bm v,t\right)d\bm x d \bm v \label{rho0}\\
\left< \bm u \right>(t)
&=
  \frac{1}{ \left< \rho\right> b_{r}^{3}}
  \int_{-b_{x}}^{b_{x}}
  \cdots
  \int_{-b_{v}}^{b_{v}}
  \bm v
  f \left(\bm x,\bm v,t\right)d\bm x d \bm v, \quad \txt{ and }\label{u0} \\
\left< T\right>(t)
&=
  \frac{1}{\left< \rho\right>  b_{x}^{3} R_{S}}
  \int_{-b_{x}}^{b_{x}}
  \cdots
  \int_{-b_{v}}^{b_{v}}\left\|\bm v - \left< \bm u \right>\right\|_2^2
  f \left(\bm x,\bm v,t\right)d\bm x d \bm v\label{T0}
\end{align}
In Figure~\ref{fig:nUTConst} we plot\eqref{rho0}-\eqref{T0} versus 
time. It is seen that such quantities are indeed constants, i.e., 
ALS iterations preserve the average density, velocity and temperature.

\subsubsection{Relaxation to Statistical Equilibrium}
\begin{figure}[t]
\center
\includegraphics[width=0.8\textwidth]{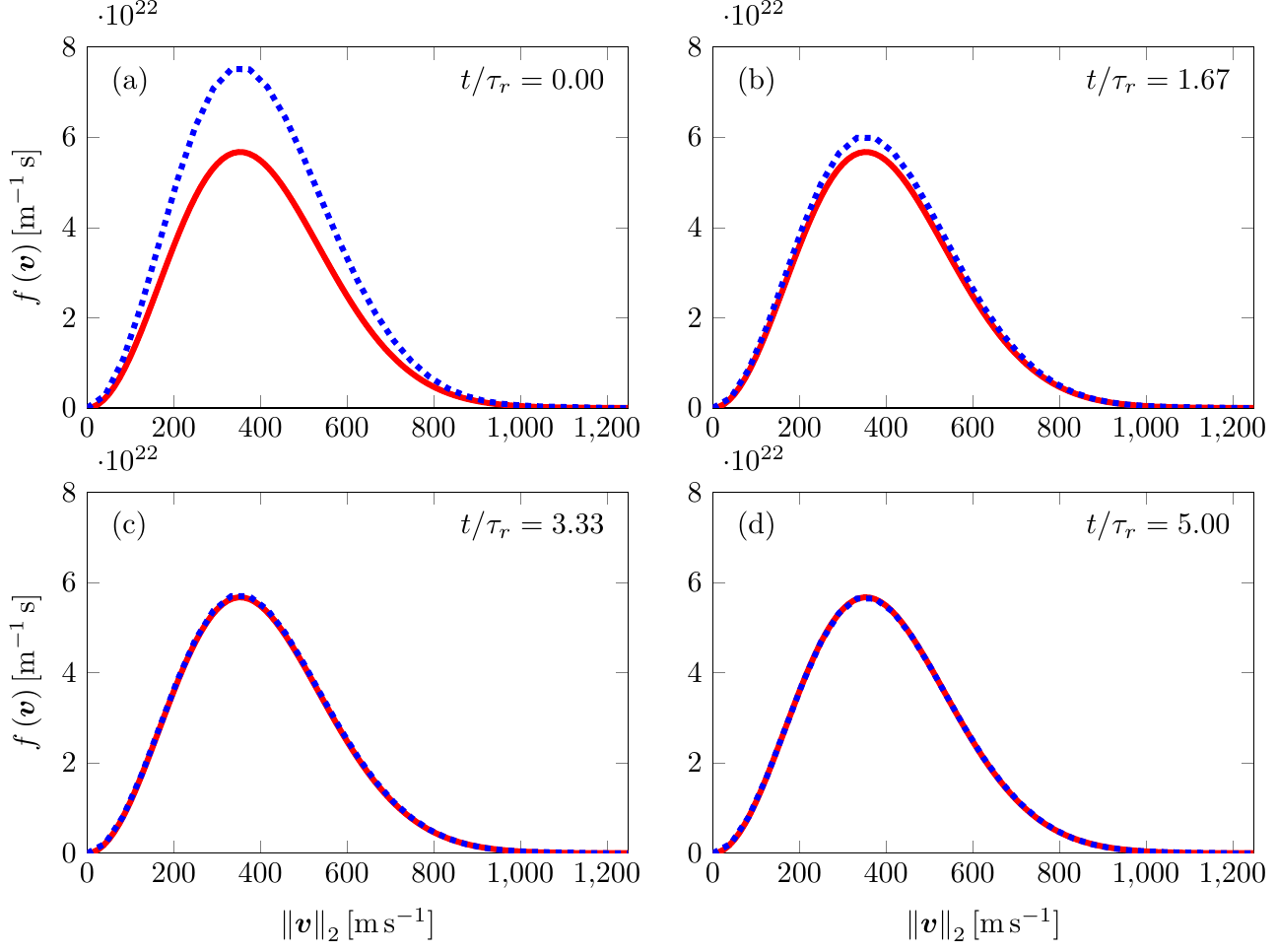}
\caption{Relaxation to statistical equilibrium in the linearized Boltzmann-BGK model \eqref{BGK_LIN}. Shown are time snapshots of the marginal PDF $p(\bm v,t)$ 
we obtain by using the parallel ALS algorithm we proposed in Section \ref{sec:BoltzmanImplicit} (\protect\includegraphics{eqNum}). 
Note that the initial condition $f(\bm v,0)$ is a slight perturbation of the 
Maxwell-Boltzmann distribution (\protect\includegraphics{eqAna}). 
Note that $f(\bm v,t)$ converges to $f_{eq}(\bm v)$ 
exponentially fast in time (see Figure \ref{fig:NMAEDecay}), 
consistently with results of perturbation analysis.}
\label{fig:cDecay}
\end{figure}

\begin{figure}[t]
\centering
\subfloat{
\includegraphics[width=0.30\textwidth]{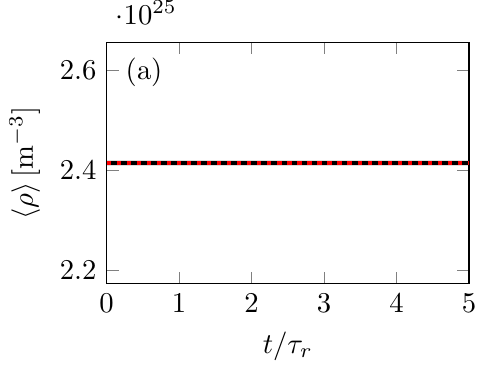}
}
\subfloat{
\includegraphics[width=0.30\textwidth]{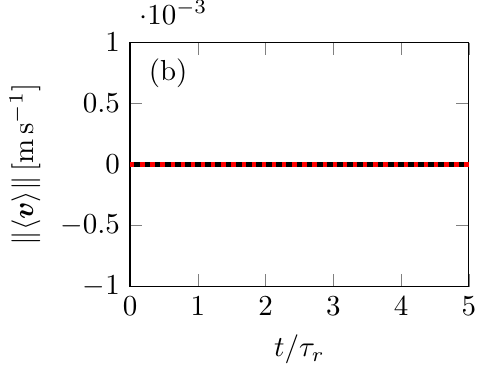}
}
\subfloat{
\includegraphics[width=0.30\textwidth]{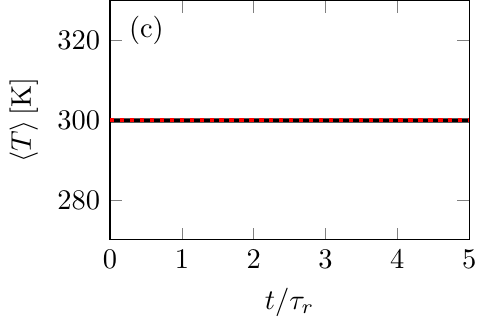}
} \\
\caption{Conservation of the average density $\left<\rho\right>$,  absolute average velocity $\left\|\left<\bm v\right>\right\|_2$, and average temperature 
$\left< T\right>$ during relaxation to statistical equilibrium. Shown are results 
obtained by using the parallel ALS algorithm we developed in Section \ref{sec:BoltzmanImplicit} (\protect\includegraphics{num}), 
and the equilibrium initial values (\protect\includegraphics{eql}).}
\label{fig:nUTDecay}
\end{figure}

\begin{figure}[t]
\centering
\includegraphics[width=\textwidth]{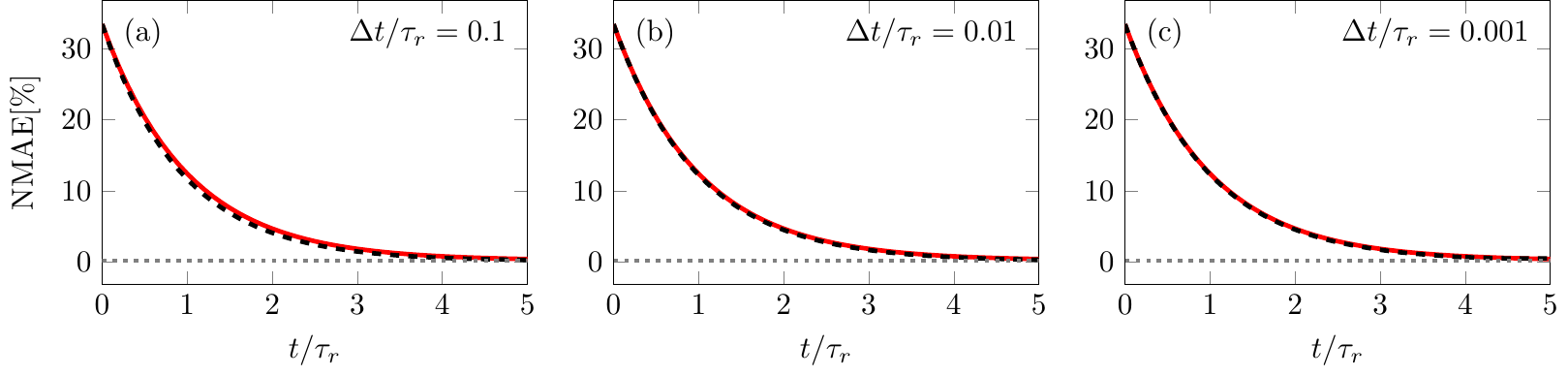}
\caption{Relaxation to statistical equilibrium. 
Normalized Mean Absolute Error (NMAE) between the
Maxwell-Boltzmann distribution $f_{eq}(\bm v)$ 
and the numerical solution we obtain with the parallel ALS 
algorithm we proposed in Section \ref{sec:BoltzmanImplicit}
(\protect\includegraphics{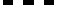}). 
We also plot the theoretically predicted exponential decay 
(\protect\includegraphics{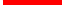}), and the best 
approximation error that we obtained based on 
approximating the equilibrium distribution $f_{eq}$ with 
a canonical tensor decomposition \eqref{functional-SSE} 
with $Q=11$ modes (\protect\includegraphics{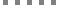}), 
$Q=11$ being the number of modes employed in 
the transient simulation.}
\label{fig:NMAEDecay}
\end{figure}

In this section we consider relaxation statistical equilibrium in the 
linearized Boltzmann-BGK model \eqref{BGK_LIN}. This  
allows us to study the accuracy and computational efficiency of the 
proposed ALS algorithm in transient dynamics. To this end, 
we consider the initial condition 
\begin{equation}
f\left(\bm x,\bm v,0\right)
= f_{eq} \left(1 + \epsilon \prod_{k=1}^{3} \cos\left(2 \pi \frac{x_{k}}{b_{x}}\right)\right) 
\label{ICf}
\end{equation}
with $\epsilon = 0.3$. Such initial condition is a slight 
perturbation of the Maxwell-Boltzmann equilibrium distribution 
(see Figure \ref{fig:cDecay}(a)). The simulation results we obtain 
are shown in Figures \ref{fig:cDecay}(b)-(d), where we plot 
the time evolution of the marginal PDF $f(\bm v,t)$ versus the 
dimensionless time $t/\tau_R$. It is seen that 
the initial condition \eqref{ICf} is in the basin of
attraction of $f_{eq}$, i.e., the ALS algorithm sends 
$f\left(\bm x,\bm v,0\right)$ into $f_{eq}\left(\bm v\right)$ after 
a transient. Moreover, $f(\bm v,t)$ is attracted by $f_{eq}(\bm v)$ 
exponentially fast in time at a rate $\tau_{R} = 1/\nu$ 
(see Figure \ref{fig:NMAEDecay}), consistently with 
results of perturbation analysis.
In Figure \ref{fig:nUTDecay} we demonstrate that the parallel ALS algorithm we developed conserves the average density of particle, velocity (momentum) and temperature throughout the transient that yields statistical equilibrium.   

In Figure \ref{fig:NMAEDecay} we plot the the Normalized Mean Absolute Error
(NMAE) between the Maxwell-Boltzmann distribution $f_{eq}(\bm v)$ 
and the distribution we computed numerically by using the proposed 
ALS algorithm. Specifically we plot NMAE versus the dimensionless 
time $t/\tau_{R}$ for different time-steps $\Delta t$ we employed
in our simulations.  
The goal of such study is to assess the sensitivity 
of parallel ALS iterations on the time step $\Delta t$, and in turn on 
the total number of steps for a fixed integration time $T=5\tau_R$. 
Results of Figure \ref{fig:NMAEDecay} demonstrate that the decay of the 
distribution function $f(\bm x,\bm v,t)$ to the its equilibrium state 
is indeed exponential. Moreover, we see that 
the ALS algorithm is robust to $\Delta t$, which implies that 
the such algorithm is not sensitive to the total number of time steps 
we consider within a fixed integration time. 

\begin{figure}
\centering
\subfloat{
\includegraphics[width=0.5\textwidth]{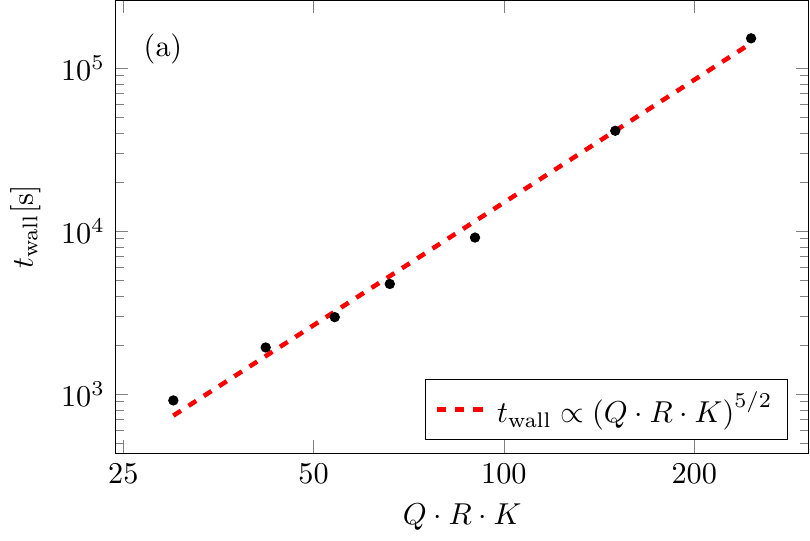}
}
\subfloat{
\includegraphics[width=0.5\textwidth]{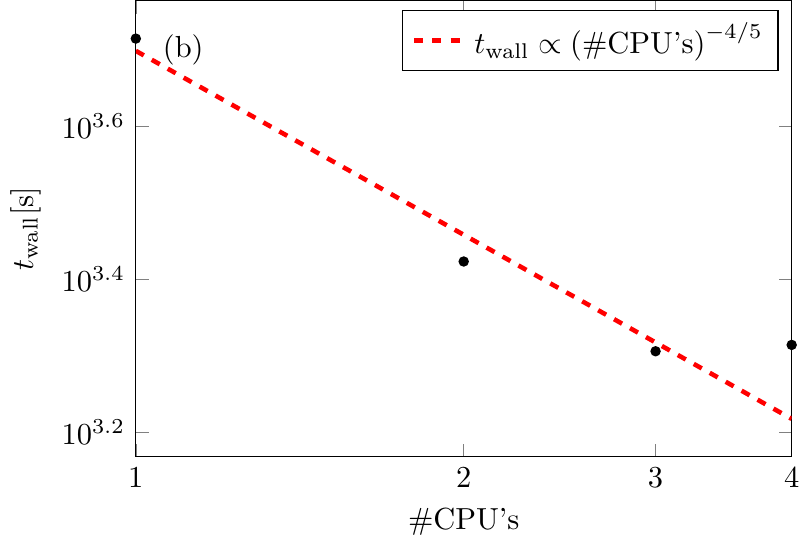}
}
\caption{Wall time $t_{\txt{wall}}$ versus the number of degrees 
of freedom of the ALS simulation (a), and versus the number 
CPUs (b).}
\label{fig:tWall}
\end{figure}

Finally, we'd like to address scalability of the parallel ALS algorithm 
we developed in  \eqref{app:parallelization}., i.e., 
performance relative to the number of degrees of freedom 
and number of CPUs. 
To this end, we analyze here a few test cases we ran on a small workstation 
with 4 CPUs. The initial condition is set as $f_{eq}$ in all cases 
and we integrated the Boltzmann-BGK model \eqref{BGK_LIN} for 
$1000$ steps, with $\Delta t= 0.01 \tau_R$. The separation rank of the 
solution PDF is set to $r=1$ (see equation \eqref{functional-SSE}). 
In Figure~\ref{fig:tWall} (a) we plot the wall time $t_{\txt{wall}}$
versus the total number of degrees of freedom of the ALS algorithm, 
i.e. $6Qr$. 
The plot shows that a typical simulation of the Boltzmann-BGK model 
in 6 dimensions with $r = 1$ and $Q=11$ takes about
$t_{\txt{wall}} \approx 80\si{\min}$ on a single core to complete 
a simulation up $t = 10\, \si{\s}$. 
With such value of $Q$ there is an error of only $0.5\%$ compared to the
analytical Maxwell-Boltzmann distribution function. 
The total wall time scales with the power $5/2$ relative to 
number of degrees of freedom $6Qr$.
In Figure~\ref{fig:tWall} we show that the total wall time scales we obtain with the proposed parallel ALS algorithm scales almost linearly 
with the number of CPUs (power $-4/5$). The leveling
off of performance for $4$ cores can be attributed to 
the fact that the desktop computer we employed for our simulations 
is optimized for inter-core communication. We expect better
scaling performance of our ALS algorithm in computer clusters with 
infinity band inter-core communication.

\section{Summary}

In this paper we presented new parallel algorithms 
to solve high-dimensional partial differential equations
based on numerical tensor methods.
In particular, we studied canonical and hierarchical tensor  
expansions combined with alternating least squares
and high-order singular value decomposition. Both implicit 
and explicit time integration methods were presented 
and discussed. We demonstrated the accuracy and 
computational efficiency of the proposed methods 
in simulating the transient dynamics generated by
a linear advection equation in six dimensions plus time, 
and the well-known Boltzmann-BGK equation.
We found that the algorithms we developed are extremely fast and
accurate, even in a naive Matlab implementation. Unlike direct
simulation Monte Carlo (DSMC), the high accuracy of the numerical 
tensor methods we propose makes it suitable for studying 
transient problems. 
On the other hand, given the nature of the numerical tensor discretization, 
implementation of the proposed algorithms to complex domains is not 
straightforward.

\section*{Acknowledgements}
We would like to thank N. Urien for useful discussions. D. Venturi 
was supported by the AFOSR grant FA9550-16-1-0092.

\appendix
\renewcommand*{\appendixname}{}

\section{Boltzmann Equation and its Approximations}
\label{sec:Boltzmann}

In the classical kinetic theory of rarefied gas dynamics, 
the gas flow is described in terms of a non-negative density 
function $f(\bm x,\bm v,t)$ which 
provides the number of gas particles at time $t$ 
with velocity $\bm v\in\mathbb{R}^3$ at position  
$\bm x \in \mathbb{R}^3$. 
The density function $f$ satisfies the Boltzmann 
equation \cite{cercignani1997,villani2002,dimarco2018}.
In the absence of external forces, such 
equation can be  written as
\begin{equation}
\frac{\partial f}{\partial t}+\bm v \cdot \nabla_{\bm x} f = Q(f,f), 
\label{BLZ}
\end{equation} 
where $Q(f,f)$ is the collision integral describing 
the effects of internal forces due to particle interactions. 
From a mathematical viewpoint the collision integral is a 
functional of the particle density. Its form depends on the 
microscopic dynamics. For example, for classical  
rarefied gas flows~\cite{cercignani1997,cercignani1994},
\begin{equation}
Q(f,f)(\bm x,\bm v,t)=\int_{\mathbb{R}^3}\int_{\mathbb{S}^2}
B(\bm v,\bm v_1,\bm \omega) \left| f(\bm x,\bm v',t)f(\bm x,\bm v_1',t)-
f(\bm x,\bm v,t)f(\bm x,\bm v_1,t) \right| \text d \bm \omega \text d\bm v_1.
\label{collisionOP}
\end{equation}
In this expression, 
\begin{equation}
\bm v' = \frac{1}{2}\left(\bm v+\bm v_1+\left\|\bm v- \bm v_1 \right\|_2\bm w\right),\qquad 
\bm v_1' = \frac{1}{2}\left(\bm v+\bm v_1-\left\|\bm v- \bm v_1 \right\|_2\bm w\right)
\end{equation}
and $\{\bm v,\bm v_1\}$ represent, respectively, 
the velocities of two particles before and after the 
collision, $\bm \omega$ is a vector on the unit sphere $\mathbb{S}^2$, 
and  $B(\bm v,\bm v_1,\bm \omega)$ is the collision kernel. 
Such kernel is a non-negative function depending on the (Euclidean) 
norm $\left\|v-v_1\right\|_2$ and on the scattering angle $\theta$ 
between the relative velocities before and after the collision
\begin{equation}
\cos(\theta) = \frac{(\bm v- \bm v_1)\cdot \bm \omega}{\left\|\bm v -\bm v_1\right\|_2}.
\end{equation}
Thus, $B(\bm v,\bm v_1,\bm \omega)$ becomes
\begin{equation}
B(\bm v,\bm v_1,\bm \omega)=\left\| \bm v - \bm v_1 \right\|_2
\sigma\left(\left\| \bm v-\bm v_1\right\|_2,\cos(\theta)\right), 
\end{equation}
$\sigma$ being the scattering cross section function \cite{cercignani1994}.
The collision operator \eqref{collisionOP} satisfies a system of 
conservation laws in the form 
\begin{equation}
\int_{\mathbb{R}^3} Q(f,f) (\bm x,\bm v, t)\psi(\bm v) \text d \bm v=0,
\label{CL}
\end{equation}
where $\psi(\bm v)$ is either $\{1,\bm v,\left\|\bm v\right\|_2^2\}$.
This yields, respectively,  conservation of mass, momentum and energy. 
Moreover, $Q(f,f)$ satisfies the Boltzmann $H$-theorem 
\begin{equation}
\int_{\mathbb{R}^3} Q(f,f)(\bm x,\bm v, t)\log\left(f(\bm x,\bm v,t)\right) 
\text d\bm v \leq 0
\end{equation}
Such theorem implies that any equilibrium distribution function, i.e.
any function $f$ for which $Q(f,f) = 0$, has the form of a locally Maxwellian
distribution
\begin{equation}
f_\text{eq}(\bm x,\bm v,t) = \frac{\rho(\bm x,t)}{(2\pi R T(\bm x,t))^{3/2}}
\exp\left(-\frac{\left\|\bm u(\bm x,t) - \bm v\right\|_2^2}{2R T(\bm x,t)}\right),
\label{LocalMaxwell}
\end{equation}
where $R$ is the gas constant, and $\rho(\bm x,t)$, $\bm u(\bm x,t)$ and $T(\bm x,t)$ are the 
density, mean velocity and temperature of the gas, defined as 
\begin{align}
\rho(\bm x,t) = &\int_{\mathbb{R}^3} f(\bm x,\bm v,t)d \bm v, \label{rho}\\ 
\bm u (\bm x,t) = &\frac{1}{\rho(\bm x,t)}\int_{\mathbb{R}^3}\bm v f(\bm x,\bm v,t) d \bm v,\label{u}\\ 
T(\bm x,t) = &\frac{1}{3\rho(\bm x,t)}\int_{\mathbb{R}^3} \left\|\bm u(\bm x,t) - \bm v \right\|_2^2 f(\bm x,\bm v,t) d \bm v \label{T}.
\end{align}
From a mathematical viewpoint, the Boltzmann equation \eqref{BLZ} is 
a nonlinear integro-differential equation for a positive scalar field in 
six dimensions plus time. By taking suitable averages over small volumes 
in position space, it can be shown that the Boltzmann equation 
is consistent with the compressible Euler's 
equations \cite{Nishida1978,caflisch1980}, and the 
Navier-Stokes equations \cite{Struchtrup2005,Levermore1996,mueller1998}.

\paragraph{BGK Approximation}
The simplest collision operator satisfying the conservation law \eqref{CL}  
and the condition 
\begin{equation}
Q(f,f)=0 \Leftrightarrow f(\bm x,\bm v,t)=f_\text{eq}(\bm x,\bm v,t)
\end{equation}
(see equation \eqref{LocalMaxwell}) was proposed 
by Bhatnagar, Gross and Krook in \cite{bhatnagar1954}. The 
corresponding model, which is known as the BGK model, is defined 
by the linear relaxation operator
\begin{equation}
Q(f,f)=\nu(\bm x,t) \left[f_\text{eq}(\bm x,\bm v,t)-f(\bm x,\bm v,t)\right] \qquad \nu(\bm x,t)>0.
\label{BoBGK}
\end{equation}
The field $\nu$ is usually set to be proportional to the density 
and the temperature of the gas \cite{Mieussens2000}
\begin{equation}
\nu(\bm x,t)=\kappa \rho(\bm x,t)T(\bm x,t)^{1-\mu}.
\label{collfreq}
\end{equation}
It can be shown that the Boltzmann-BGK model \eqref{BoBGK} 
converges to the correct Euler equations of incompressible fluid dynamics 
with the scaling $\bm x'=\epsilon \bm x$, $t'=\epsilon t$, and in the 
limit $\epsilon\rightarrow 0$.  
However, the model does not converge to the correct 
Navier-Stokes limit. The main reason is that it predicts an 
unphysical Prandtl number \cite{nassios2013}, which is 
larger that one obtained with the full collision 
operator \eqref{collisionOP}. 
The correct Navier-Stokes limit can be recovered by more 
sophisticated BGK-type models, e.g., the 
ellipsoidal statistical BGK model \cite{Andries2000}.

Following \citep{borghini2016,barichello2003, pekeris1957, loyalka1967}, 
we introduce additional simplifications, namely we assume 
that  $\nu(\bm x,t) =\nu$ is constant and assume the equilibrium density, temperature and velocity 
to be homogeneous across the spatial domain. Assuming $\bm u(\bm x,t)=\bm 0$, this yields an equilibrium distribution in~\eqref{EqMaxwell}. The last assumption effectively decouples the BGK collision 
operator from the probability density function $f(\bm x,\bm v,t)$. 
This, in turn, makes the BGK approximation linear, i.e., a six-dimensional PDE~\eqref{BGK_LIN}.
Using both the isothermal assumption $T(\bm x,t)=T$ \cite{bhatnagar1954}, and
constant density assumption $\rho(\bm x,t)=\rho$ \cite{barichello2003,
pekeris1957, loyalka1967} means that our code can only be used to obtain
solutions which are small perturbations from the equilibrium distribution.
The term $\nu \left[f_\text{eq}(\bm v)-f\right]$ in \eqref{BGK_LIN} 
represents the relaxation time approximation of the Boltzmann 
equation, $\nu$ is the collision frequency, here assumed to be 
homogeneous.

\section{Parallelization of the ALS Algorithm for the Boltzmann-BGK Equation}
\label{app:parallelization}
In this appendix we describe the parallel alternating least squares 
algorithms we developed to solve the linearized Boltzmann-BGK equation
\eqref{BGK_LIN}.
Algorithm \ref{alg:als} shows the main Alternating 
Least Square (ALS) routine while algorithm
\ref{alg:CreateMatrixM}, \ref{alg:CreateVectorGamma}, \ref{alg:ComputeBeta}, and
\ref{alg:ComputeConvergence} show the different subroutines. Each time step,
$n$, the value of $\beta_{\txt{Old}}$ is updated and then $\beta_{\txt{New}}$ is
determined in an iterative manner. Both the computations of the matrices $M_{k}$
and the vectors $\gamma_{k}$, as well as the calculation of the vectors
$\beta_{\txt{New},k}$ can be performed independently of each other, and thus can
be easily parallelized. The current implementation assigns one core per dimension,
but further optimization could be possible by implementing a parallel method to
determine $\beta_{\txt{New}}$ \cite{huang2013}.
\begin{algorithm}[t]
\caption{Parallel Alternating Least Squares (ALS) algorithm}
\label{alg:als}
\begin{algorithmic}[1]
\Procedure{Main}{}
\State{\Call{Initialization}{}} \Comment{Load variables and allocate matrices}
\State{$\beta_{\txt{New}} = \beta_{0}$} \Comment{Set initial condition}
\For{$N \gets 1:N_{\txt{Iter}}$}
\State{$\epsilon_{\txt{Conv}} \gets \epsilon_{\txt{Tol}} + 1$} \Comment{Reset convergence condition}
\State{$\beta_{\txt{Old}} \gets \beta_{\txt{New}}$}
\While{$\epsilon_{\txt{Conv}} > \epsilon_{\txt{Tol}}$}
\State{$\beta_{\txt{Int}} \gets \beta_{\txt{New}}$}
\ParFor{$k \gets 1:K$} \Comment{This for loop is executed in parallel}
\State{$\hphantom{M_{k}}\mathllap{M_{k}}      \gets \Call{CreateMatrixM    }{\beta_{\txt{New},k}                    ,k}$}
\State{$\hphantom{M_{k}}\mathllap{\gamma_{k}} \gets \Call{CreateVectorGamma}{\beta_{\txt{New},k},\beta_{\txt{Old},k},k}$}
\EndParFor
\ParFor{$k \gets 1:K$} \Comment{This for loop is executed in parallel}
\State{$\beta_{\txt{New},k} \gets \Call{ComputeBeta}{M_{k},\beta_{\txt{New},k},\gamma_{k}}$} \Comment{Compute updated $\beta_{\txt{New}}$}
\EndParFor
\State{$\epsilon_{\txt{Conv}} \gets \Call{ComputeConvergence}{\beta_{\txt{Int}},\beta_{\txt{New}}}$} \Comment{Update convergence condition}
\State{$\beta_{\txt{New}} \gets \beta_{\txt{Int}} + \brc{\beta_{\txt{New}} - \beta_{\txt{Int}}}/4$} \Comment{Update $\beta_{\txt{New}}$}
\EndWhile
\EndFor
\EndProcedure
\end{algorithmic}
\end{algorithm}
At each time step the CreateMatrixM subroutine
(algorithm~\ref{alg:CreateMatrixM}) updates the different elements of matrix
$M_{k}$ for every iteration until convergence has been reached. The algorithm
iterates over every element of the matrix and performs the multiplication and
summation of the elements of $\hat{M}$. All the integrals were performed
analytically and the results are stored in the map $\txt{M}\brc{k,l,z,s,q}$. 
Using this approach delivers significant savings in the 
computational cost of updating the matrix every iteration. 
\begin{algorithm}[t]
\caption{ALS Subroutines: CreateMatrixM}
\label{alg:CreateMatrixM}
\begin{algorithmic}[1]
\Function{CreateMatrixM}{$\beta_{\txt{New},k},k$}
%
\For{$l \gets 1:R_{G}$} \Comment{Loop over all elements of $M_{k}$}
\For{$z \gets 1:R_{G}$}
\For{$s \gets -\abs{Q}/2:\abs{Q}/2$}
\For{$q \gets -\abs{Q}/2:\abs{Q}/2$}
\For{$k' \gets 1:K$} \Comment{Loop over all elements of $\hat{M}$ for each element of $M_{k}$}
\For{$\eta \gets 1:r_{A}$}
\For{$\upsilon \gets 1:r_{A}$}
\State{$j \gets \brc{\eta-1} r_{A} + \upsilon$}
%
\If{k' = k}
\State{$\hat{M}_{k',j} \gets \Call{M}{k',\eta,\upsilon,s,q}$}
\Else
\State{$\hat{M}_{k',j} \gets 0$}
\For{$s' \gets -\abs{Q}/2:\abs{Q}/2$}
\For{$q' \gets -\abs{Q}/2:\abs{Q}/2$}
\State{$\hat{M}_{k',j} \gets \hat{M}_{k',j} + \beta_{\txt{New},k',\brc{l-1}Q+s'} \beta_{\txt{New},k',\brc{z-1}Q+q'} \Call{M}{k',\eta,\upsilon,s',q'}$}
\EndFor
\EndFor
\EndIf
\EndFor
\EndFor
\EndFor
\State{$M_{k,\brc{z-1}Q+q,\brc{l-1}Q+s} \gets \sum_{j=1}^{r_{A}^{2}} \prod_{k'=1}^{K} \hat{M}_{k',j}$} \Comment{Fill matrix $M_{k}$}
\EndFor
\EndFor
\EndFor
\EndFor
\State{\Return{$M_{k}$}}
\EndFunction
\end{algorithmic}
\end{algorithm}
The subroutine \texttt{CreateVectorGamma} (algorithm~\ref{alg:CreateVectorGamma}) updates
the vector $\gamma_{k}$. Similarly to the subroutine for the creation of matrix
$M_{k}$ this algorithm iterates over all elements of the vector, and fills each
of them using the values of $\beta_{\txt{New}}$ and $\beta_{\txt{Old}}$, and
pre-calculated maps. 


\begin{algorithm}[t]
\caption{ALS Subroutines: CreateVectorGamma}
\label{alg:CreateVectorGamma}
\begin{algorithmic}[1]
\Function{CreateVectorGamma}{$\beta_{\txt{New},k},\beta_{\txt{Old},k},k$}
\For{$z \gets 1:R_{G}$} \Comment{Loop over all elements of $\gamma_{k}$}
\For{$q \gets -\abs{Q}/2:\abs{Q}/2$}
\For{$k' \gets 1:K$} \Comment{Loop over all elements of $\hat{N}$ for each element of $\gamma_{k}$}
\For{$l' \gets 1:R_{G}$}
\For{$\eta \gets 1:r_{A}$}
\For{$\upsilon \gets 1:r_{A}$}
\State{$j \gets \brc{\eta-1} r_{A} + \upsilon$}
%
\If{k' = k}
\State{$\hat{N}_{k',l',j} \gets 0$}
\For{$s' \gets -\abs{Q}/2:\abs{Q}/2Q$}
\State{$\hat{N}_{k',l',j} \gets \hat{N}_{k',l',j} + \beta_{\txt{Old},k',\brc{l'-1}Q+s'} \Call{N}{k',\eta,\upsilon,s',q}$}
\EndFor
\Else
\State{$\hat{N}_{k',l',j} \gets 0$}
\For{$s' \gets -\abs{Q}/2:\abs{Q}/2Q$}
\For{$q' \gets -\abs{Q}/2:\abs{Q}/2$}
\State{$\hat{N}_{k',l',j} \gets \hat{N}_{k',l',j} + \beta_{\txt{Old},k',\brc{l'-1}Q+s'} \beta_{\txt{New},k',\brc{z-1}Q+q'} \Call{N}{k',\eta,\upsilon,s',q'}$}
\EndFor
\EndFor
\EndIf
\EndFor
\EndFor
\EndFor
\EndFor
\State{$N_{k,\brc{z-1}Q+q} \gets \sum_{l'=1}^{R_{G}} \sum_{j=1}^{r_{A}^{2}} \prod_{k'=1}^{K} \hat{N}_{k',l',j}$} \Comment{Fill matrix $N_{k}$}
\For{$k' \gets 1:K$} \Comment{Loop over all elements of $\hat{O}$ for each element of $\gamma_{k}$}
\For{$\eta \gets 1:r_{A}$}
\If{k' = k}
\State{$\hat{O}_{k',\eta} \gets \Call{O}{k',\eta,q}$}
\Else
\State{$\hat{O}_{k',\eta} \gets 0$}
\For{$q' \gets -\abs{Q}/2:\abs{Q}/2$}
\State{$\hat{O}_{k',\eta} \gets \hat{O}_{k',\eta} + \beta_{\txt{New},k',\brc{z-1}Q+q'} \Call{O}{k',\eta,q'}$}
\EndFor
\EndIf
\EndFor
\EndFor
\State{$O_{k,\brc{z-1}Q+q} \gets \sum_{\eta=1}^{r_{A}} \prod_{k'=1}^{K} \hat{O}_{k',\eta}$} \Comment{Fill matrix $O_{k}$}
\EndFor
\EndFor
\State{$\gamma_{k} = N_{k}+\Delta t \; \nu \; O_{k}$}
\State{\Return{$\gamma_{k}$}}
\EndFunction
\end{algorithmic}
\end{algorithm}
With the completion of the calculation of $M_{k}$ and $\gamma_{k}$,
$\beta_{\txt{New}}$ can be updated using the \texttt{ComputeBeta} subroutine
(algorithm~\ref{alg:ComputeBeta}). Matrix $M_{k}$, depending on the rank
$r_{G}$, and series truncation $Q$, can become very large. In addition, it is
not known whether the matrix is invertible. This is why the system $M_{k}
\beta_{k} = \gamma_{k}$ is solved using the iterative least square method, LSQR
\cite{paige1982,barrett1994}. For every iteration the value for
$\beta_{\txt{New},k}$ which was determined in the previous iteration is used as
the initial guess for the algorithm.

\begin{algorithm}[t]
\caption{ALS Subroutines: ComputeBeta}
\label{alg:ComputeBeta}
\begin{algorithmic}[1]
\Function{ComputeBeta}{$M_{k},\beta_{\txt{New},k},\gamma_{k}$}
\State{$\beta_{\txt{New},k} \gets
\Call{Lsqr}{M_{k},\beta_{\txt{New},k},\gamma_{k}}$} \Comment{The system
$\gamma_{k} - M_{k} \beta_{\txt{New},k} = \txt{Res}$ is solved using the least
squares method. $\beta_{\txt{New},k}$ is used as an initial guess before being updated} 
\State{\Return{$\beta_{\txt{New},k}$}}
\EndFunction
\end{algorithmic}
\end{algorithm}

Convergence of the system is calculated using the \texttt{ComputeConvergence} subroutine
(algorithm~\ref{alg:ComputeConvergence}). Because it is very computationally
expensive to calculate the full residual every iteration, instead the
convergence is determined using:
\begin{equation}
  \epsilon_{\txt{Conv}} 
= 
  \max{\brc{  \frac{\nor{\beta_{\txt{New},k} -
  \beta_{\txt{Int},k}}}{\nor{\beta_{\txt{Int},k}}} }}
\end{equation}
and when $\epsilon_{\txt{Conv}} < \epsilon_{\txt{Tol}}$, the while loop in
algorithm~\ref{alg:als} is allowed to exit. In its currently implementation the
algorithm is able to reach convergence using a fixed Rank G, $r_{G} = 2$.
However, a future implementation could implement a dynamic Rank G for improved
accuracy. Every iteration the new value of
$\beta_{\txt{New}}$ is calculated according to:
\begin{equation}
  \beta_{\txt{New}} 
\gets
  \beta_{\txt{Int}} + \frac{{\beta_{\txt{New}} - \beta_{\txt{Int}}}}{\Delta
  \beta}
\end{equation}
with $\Delta \beta = 4$. To find the location where the residual reaches a
minimum, it is important to gradually approach this location and not update
$\beta_{\txt{New}}$ too aggressively by taking a smaller value of $\Delta \beta$.
This causes overshooting of the minimum and the prevents convergence of
$\beta_{\txt{New}}$ and thus minimization of  the residual.

\begin{algorithm}[t]
\caption{ALS Subroutines: ComputeConvergence}
\label{alg:ComputeConvergence}
\begin{algorithmic}[1]
\Function{ComputeConvergence}{$\beta_{\txt{Int}},\beta_{\txt{New}}$}
\For{$k \gets 1:K$}
  \State{$\beta_{\txt{Norm},k} \gets {\nor{\beta_{\txt{New},k} - \beta_{\txt{Int},k}}}/{\nor{\beta_{\txt{Int},k}}}$}
\EndFor
\State{ \Return{$\max{\brc{\beta_{\txt{Norm}}}}$}}
\EndFunction
\end{algorithmic}
\end{algorithm}

\end{document}